\documentclass[11pt]{article}
\usepackage{array,graphicx, xcolor, lipsum, bibentry}
\usepackage[margin=3cm]{geometry}
\usepackage{tikz}
\usetikzlibrary{shapes.callouts}
\usetikzlibrary{calc}
\usepackage{amsmath}
\usepackage{amsfonts}
\usepackage{eurosym}
\usepackage{hyperref}
\usepackage[authoryear,round]{natbib}
\usepackage{cleveref}
\usepackage{url}
\usepackage{subcaption,caption}
\usepackage{multirow}
\usepackage[title,titletoc]{appendix}
\usepackage{longtable}
\usepackage{authblk}
\usepackage{booktabs}
\providecommand{\keywords}[1]{\textbf{\textit{Keywords:}} #1}
\crefname{appsec}{Appendix}{Appendices}

\newcommand{\mc}[1] {\mathcal{#1}}

\allowdisplaybreaks

\title{Maximizing performance with an eye on the finances: a chance-constrained model for football transfer market decisions}

\author[1]{G. Pantuso}

\author[2]{L.M. Hvattum}

\affil[1]{\small Department of Mathematical Sciences, University of Copenhagen\\
Universitetsparken~5, 2100 Copenhagen, Denmark\\
email: \href{mailto:gp@math.ku.dk}{gp@math.ku.dk}\\
tel:~+4535335608. Corresponding author.
}

\affil[2]{Faculty of Logistics, Molde University College\\
Britvegen 2, 6410 Molde, Norway\\
email: \href{mailto:hvattum@himolde.no}{hvattum@himolde.no}\\
tel:~+4771214223
}

\date{}

\begin{document}
\maketitle
\begin{abstract}
  Composing a team of professional players is among the most crucial decisions in association football. Nevertheless, transfer market decisions are often based on myopic objectives and are questionable from a financial point of view. This paper introduces a chance-constrained model to provide analytic support to club managers during transfer windows. The model seeks a top-performing team while adapting to different budgets and financial-risk profiles. In addition, it provides a new rating system that is able to numerically reflect the on-field performance of football players and thus contribute to an objective assessment of football players. The model and rating system are tested on a case study based on real market data. The data from the case study are available online for the benefit of future research. 
\end{abstract}
\keywords{association football, team composition, sport, sports management, stochastic programming, integer programming}

\section{Introduction}\label{sec:Intro}
Composing a team of players is among the most crucial decisions a football club's manager is required to make.
In fact, the main component of a football club's costs is expenditure on players, through both wages and transfer fees \citep{DobG01}.
However, the analysis performed by \citet{KupS18} illustrates that transfer decisions are often based on myopic objectives, impulsive reactions, and are overly influenced by factors such as recent performances \citep{Lew04}, particularly in big tournaments, nationality, and even hair color. With few exceptions, football clubs are in general bad businesses, rarely making profits and most often accumulating considerable debts. The authors explain that poor management of football clubs is tolerated because, unlike in most industries, clubs in practice never go bankrupt. History shows that creditors rarely claim their credits to the extent of causing the default of the club. In other words: \textit{``no bank manager or tax collector wants to say: The century-old local club is closing. I am turning off the lights''} \citep{KupS18}. In practice, clubs eventually find someone who bails them out, change management, narrow the budget, and forces them to restart by competing at lower standards (e.g., being relegated).

According to UEFA, this dramatic scenario is changing for the better. In their \textit{Club Licensing Benchmarking Report} for the financial year 2016 \citep{UEFA18}, UEFA reports that the combined debt of Europe's top-division clubs (which includes net loaning and net player transfer balance) has decreased from 60\% of revenue in 2010 to 35\% of revenue in 2016, partially as a consequence of the introduction of UEFA Financial Fair Play rules. In addition, the report shows that profits are also increasing, generating resources to be reinvested in the football business, e.g., to reinforce the team. Furthermore, a number of football clubs are listed in stock exchanges (see, \citet{KPMG17}), with the notable case of Juventus FC entering the FTSE-MIB index which tracks the performance of the 40 main shares in the Italian stock market \citep{Sole18,Van18}. In this paper, we also advocate for the stability and good management of football clubs. Even if traditionally poorly managed clubs have often found a benefactor, more solid decision making would prevent financial distress, improve the capacity to generate resources to reinvest in the team, and eventually spare worries to the club's supporters. In particular, in this paper we introduce analytic methods for supporting transfer decisions.

In the research literature, \citet{Pan17} proposed a stochastic programming model with the scope of maximizing the expected market value of the team. The findings of the study confirm analytically the recipe provided by \citet{KupS18}: a steady growth in the team value is associated with fewer transfers, timely selling old players (they are often overrated), and investments in your prospects. However, \citet{Pan17} did not explicitly account for the players' performances. While the strategy of maximizing the expected market value of the team might benefit the club in the long-run, it may as well contrast with the short-term requirement of meeting competitive goals. In fact, most managers are often evaluated by matches won. \citet{PaZh19} presented a deterministic integer programming model, maximizing a weighted sum of player values, adjusted for age and rating, minus net transfer expenses. Their model covered several time periods, and was tested on twelve teams from the English Premier League.

The fact that not all football clubs act as profit maximizers has been mentioned in scientific literature on many occasions, starting with \citet{Sl71} who suggested that European clubs behave as utility maximizers, with a utility function that contains other variables in addition to profit. \citet{Ke96} introduced win maximization, with the consequence that clubs should hire the best players within the limits of their budget. \citet{Ra97} considered clubs to maximize a linear combination of wins and the profit level, with different clubs having a different weight to balance the two criteria. \citet{Ke06} concluded that most clubs are interested in more than making profit, but also that they do not want to win at any cost. This in turn translates into the requirement of hiring top-performers who can immediately contribute to on-field successes, but at the same time keeping an eye on the financial performance of the club.

Measuring on-field performance means associating a numerical value to the contribution given by the player to the team. This can been done using different methods \citep{Sz15}, one of which is the plus-minus rating. These initially consisted of recording the goals scored minus the goals conceded from the perspective of each player, and were applied to ice hockey players. \citet{winston} showed how the principle could be applied in basketball by calculating adjusted plus-minus ratings, which are determined by multiple linear regression. This allows the ratings to compensate for the teammates and opponents of each player. The next important development was to use ridge regression instead of the method of ordinary least squares to estimate the regression model, as proposed by \citet{Si10} for basketball, and later \citet{macdonald12} for ice hockey. This is sometimes known as regularized plus-minus, and was adapted to association football by \citet{SaHv15}, with later improvements by \citet{SaHv17}. \citet{Hv19} provided an overview of the different developments made for plus-minus ratings, covering association football as well as other team sports. One of the contributions of this paper is to present an improved regularized plus-minus for association football, obtained by adding several novel features.

To account for both on-field and financial performances when composing a football team, this paper provides a chance-constrained \citep{ChaC59} mixed-integer optimization model. The objective of the model is that of finding the mix of players with the highest sum of individual player ratings. The selected mix of players must provide the skills required by the coach. In addition, the total net expenditure in transfer fees must respect the given budget. Finally, the future value of the players in the team must remain above a specified threshold with a given probability. The latter condition is enforced by a chance constraint. In contrast to the model of \citet{Pan17}, the new model does not primarily maximize the market value of the team, but rather a performance-based rating. 

The proposed model can support football managers in the course of transfer market sessions. In this phase, managers engage in negotiations with other clubs with the scope of transferring target players. Key parameters such as player values, transaction costs and salaries are typically updated as negotiations progress. The proposed model may be employed at every stage of the negotiation to assess the risk of the potential move and its impact on the remaining transfers. By using such models, football managers may obtain an analytical validation of their decisions and avoid the myopic and biased decision making which, as we commented above, has often characterized transfer market decisions. In order to use a similar model, a football club needs some sort of analytical expertise. As an example, the club must be able to compute ratings, and have an outlook on the future value of the players negotiated. While it is very likely that the majority of the professional clubs are currently not ready for such a transition, it also true that an increasing number of football clubs are starting to use some sort of analytics and data science techniques, and that vast volumes of data are being collected both on match events as well as on transactions. In this sense, the model we propose is one of the first attempts to build a decision support tool for aiding football clubs financial decisions. As football clubs will reach some maturity in the adoption of analytics techniques, the present model can be further extended to fully capture the complexity faced by football managers.

The contributions of this paper can be summarized as follows:
\begin{itemize}
\item A novel optimization problem which supports team composition decisions while accounting both for the on-field and financial performance of the club.
\item An improved player-rating system which significantly improves on state-of-the-art plus-minus ratings for football players.
\item An extensive computational study based on real transfer market data which highlights the results achievable with the new optimization model and rating system, as well as the differences between the solutions provided by our model and the solution to an existing model from the literature. Furthermore, the case studies used in the computational study are made available online at \url{https://github.com/GioPan/instancesFTCP} in order to facilitate future research on the topic.
\end{itemize}

This paper is organized as follows. In \Cref{sec:Problem} we provide a more thorough description of the problem and provide a mathematical formulation in the form of a chance-constrained mixed-integer program. In \Cref{sec:PlayerRatings} we introduce a novel plus-minus player-rating system. In \Cref{sec:casestudies} we introduce and explain the case studies. In \Cref{sec:results} we analyze the decisions obtained with our model based on historical English Premier League data, and compare the decisions to those obtained from the model provided by \citet{Pan17}. Finally, conclusions are drawn in \Cref{sec:Conclusions}.

\section{Problem Definition}\label{sec:Problem}
The manager of a football club has to decide how to invest the available budget $B$ to compose
a team of football players. Particularly, the manager's decisions include which players
to buy or loan from other clubs, and which players to sell or loan out to other clubs.
Let $\mc{P}$ be the set of all players considered, both those currently in the team and those
the club is considering buying or loaning. The latter will be referred to as \textit{target players}. Let parameter $Y_p$ be $1$ if the player $p$
belongs to the club at the planning phase, and $0$ otherwise. As decided by (inter)national football associations, each team must be composed of at least $N$ players.
The players in excess cannot participate in competitions and we therefore assume they must be sent to other teams on a loan agreement. The alternative would be that a player remains at the club but is not registered for official competitions. Though possible, this situation is undesirable.

The players in the team must cover a set of roles $\mc{R}$.
In particular, let $\underline{N}_r$ and $\overline{N}_r$ be the minimum and maximum number of players,
respectively, in role $r$. A role is, in general, a well defined set of technical and personal
characteristics of the player, such as the position on the field of play, the nationality, the speed, or strength.
The players required in a given role are typically decided by the coach when the role corresponds to
a technical characteristic. However, when the role defines a personal characteristic such as age or nationality,
national or international regulations may specify how many players with those characteristics a club may employ.
As an example, clubs competing in the Italian Serie A may not employ more than three non-EU citizens, and must employ at least
four players trained in the academy of an Italian club. Let $\mc{P}_r$ be the set of players having role $r\in\mc{R}$.
Notice that players might have more than one role so that $\bigcap_{r\in \mc{R}}\mc{P}_r\neq \emptyset$. Also, while in reality a player performs better in one role rather than another, here we consider the role a binary trait, that is, a player may occupy a given role or not, independently of how well they perform in that role. It is then a decision of the coach to employ the player in the role where they perform best.

For each target player $p$ the club is assumed to know the current purchase price $V^P_{p}$ and loan fee $V^{B}_{p}$. A target player may also be a player currently in the club's own youth team or second team who is considered for a promotion. In this case, the purchase price might be set to either zero or the opportunity cost generated by the lost sale of the player. Similarly,
for each player $p$ currently in the team, the club knows the current selling price $V^S_{p}$ and loan fee $V^{L}_{p}$. Observe that current purchase and selling prices as well as loan fees are known to the decision maker. They are either the result of a negotiation or is information which a player's agent can obtain from the owning club. In certain cases they are even explicitly stated in contracts. Notice also that today we observe several forms of payments, such as in instalments, and with bonuses conditional on the achievement of sports results. These, more involved, forms of payment may play an important role for the outcome of a negotiation, for a club's financial stability and for complying with regulations such as the UEFA Financial Fair Play. However, in our context we simply consider the actualized sum of payments since we are only concerned with ensuring, with a certain probability, that the value of investments exceeds a given threshold in the future.

However, the future market value of the player is uncertain and dependent on several unpredictable factors such as fitness,
injuries and successes. Let random variable $\tilde{V}_p$ represent the market value of player $p$ at a selected future point in time, e.g., at the beginning of next season. We assume that the probability distribution of $\tilde{V}_p$ is known to the decision maker. Such a distribution may be the result of a forecast based on historical data as done in \Cref{subsec:regression}, or more simply the outcome of expert opinions.

Let $W_p$ be the rating of player $p$, corresponding to a measure of the on-field performance of the player.
The objective of the club is that of composing a team with the highest \textit{rating}, such that
the size of the team is respected, the number of players in each role is respected, the budget is not
exceeded, and that the probability that the market value of the team exceeds a threshold $V$ is higher than $\alpha$.

Let decision variable $y_p$ be equal to $1$ if player $p$ belongs to the club (not on a loan agreement) at the end of the focal \textit{Transfer Market Window} (TMW), and $0$ otherwise. A TMW is a period of time, decided by national and international football associations, during which clubs are allowed to transfer players. Variables $y^B_{p}$ and $y^S_{p}$ are equal to $1$ if player $p$ is bought or sold during the TMW, respectively, and $0$ otherwise. Similarly, variables $x^L_{p}$ and $x^B_{p}$ are equal to $1$ if player $p$ is respectively leaving or arriving in a loan agreement during the TMW, and $0$ otherwise. The \textit{Football Team Composition Problem} (FTCP) is thus expressed by the following optimization problem.

\begin{subequations}
\label{eq:FTCP}
\begin{align}
  \label{eq:obj1}\max &\sum_{p\in \mc{P}}W_p(y_p+x^B_p-x^L_p)\\
  \label{eq:balance0}\text{s.t. }&y_{p} - y^B_{p} + y^S_{p} = Y_p  & p \in \mc{P},\\
  \label{eq:max_squad_size}&\sum_{p\in \mc{P}}(y_{p}+x^B_p-x^L_{p}) = N, & \\
  \label{eq:role_covering1}&\sum_{p\in \mc{P}_r}(y_{p}+x_p^B-x^L_p)\geq \underline{N}_r & r\in \mc{R},\\
  \label{eq:role_covering2}&\sum_{p\in \mc{P}_r}(y_{p}+x_p^B-x^L_p)\leq \overline{N}_r & r\in \mc{R},\\
  \label{eq:loan_in_if_not_owned}&x^B_p + y^B_p\leq 1-Y_{p}  & p\in \mc{P},\\
  \label{eq:loan_out_if_owned}&x^L_{p} +y^S_p\leq Y_p & p\in \mc{P},\\
  \label{eq:budget_limit}&\sum_{p \in \mc{P}} \left(V^P_{p}y^B_{p}+V^{B}_{p}x^B_{p}-V^S_{p}y^S_{p}-V^{L}_{p}x^L_{p}\right) \leq B, & \\
  \label{eq:chance}&P\bigg(\sum_{p\in \mc{P}}\tilde{V}_py_p \geq V\bigg)\geq \alpha,&\\
  \label{eq:range}&y_{p},y^B_{p},y^S_{p},x^B_{p},x^L_{p} \in \{0,1\} & p\in \mc{P}.
\end{align}
\end{subequations}

The objective function \eqref{eq:obj1} represents the performance of the team
as described by the sum of the ratings of the players competing for the team. The model maximizes the total rating of the entire team, irrespective of the lineup chosen for each specific match. That is, the model provides the coach the best possible team, and then it is the coach's duty to decide lineups for each match day.
Constraints \eqref{eq:balance0} ensure that a player belongs to the team if he
has been bought or if the player was in the team before the opening of the TMW
and has not been sold. Constraints \eqref{eq:max_squad_size} ensure that the club
registers exactly $N$ players for competitions, while constraints \eqref{eq:role_covering1}
and \eqref{eq:role_covering2} ensure that each role is covered by the necessary number of players.
Constraints \eqref{eq:loan_in_if_not_owned} and \eqref{eq:loan_out_if_owned} ensure that
a player is loaned only if not owned, and sent on a loan only if owned, respectively.
Constraints \eqref{eq:budget_limit} ensure that expenses for obtaining
players can be financed by players sold, players loaned out, or a separate budget. Constraint \eqref{eq:chance}
ensures that the probability that the future value of the team exceeds a threshold $V$
(e.g., the current value of the team) is higher than $\alpha$, with $\alpha$ reflecting
the financial risk attitude of the club. Finally, constraints \eqref{eq:range} set the
domain for the decision variables. In addition, similarly to \citet{Pan17}, it is
possible to fix the value of a subset of the variables to indicate whether a player cannot
be sold, bought, or moved on loan.

A possible limitation of model \eqref{eq:FTCP} is that it does not explicitly take into account the salary of players, which may also play an important role in the composition of a team. As an example, salaries may consume the available budget $B$. In this case, it is possible to modify constraints \eqref{eq:budget_limit} by subtracting from the budget the salaries of the players in the team (purchased or in a loan agreement), and adding the salaries of the players sold or leaving in a loan agreement for the next season. Similarly, salaries may affect the chance constraint \eqref{eq:chance}. As an example, a club may want to enforce, with a certain probability, that the future value of the players, net of the salaries paid, exceeds a certain threshold $V$. Furthermore, in principle, the model allows an uneven distribution of ratings among the team. That is, it is possible that the resulting team is composed of a few highly rated players while the remaining players have very low ratings. However, this scenario is extremely unlikely under the reasonable assumption that target lists are somewhat homogeneous in terms of rating. That is, we do not expect that a team targets both a top performer and a very poor performer. Nevertheless, extensions of the model might be considered which either distinguish between the ratings of the ideal starting line-up and the ratings of the available substitutes, or ensures a fair distribution of ratings across the entire team, or across roles, in order to foster internal competition. We leave these possible extensions to future research and, in  what follows, we employ model \eqref{eq:FTCP} without any modifications.

\section{Player Ratings}\label{sec:PlayerRatings}

Multiple linear regression models, as used to calculate adjusted plus-minus ratings, are typically stated using $y$ to denote the dependent variable, and $y_i$ being the value of the dependent variable in observation $i$. A set $\mc{V}$ of independent variables, denoted by $x_{j}$ for $j \in \mc{V}$ and with values $x_{ij}$ for observation $i$, are assumed to be related to the dependent variable such that
\begin{align}
\sum_{j \in \mc{V}} \beta_{j} x_{ij} = y_{i} + \epsilon_{i} \nonumber 
\end{align}

\noindent where $\beta_j$ are parameters describing the relationship between the independent variables and the dependent variable, and $\epsilon_i$ is an error term. When using ordinary least squares to estimate the values of $\beta_i$, one is essentially solving the following unconstrained quadratic program, with $n$ being the number of observations:
\begin{align}
\min_{\beta} \left \{ \sum_{i=1}^n ( \sum_{j \in \mc{V}} x_{ij} \beta_{j} - y_{i} )^2 \right \} \nonumber 
\end{align}

For regularized plus-minus ratings, Tikhonov regularization, also known as ridge regression, is employed instead of ordinary least squares. The main purpose of this is to avoid overfitting the model as a result of collinearity, by shrinking all regression coefficients towards zero \citep{Si10}. For example, with standard adjusted plus-minus ratings, players with few minutes played are prone to being assigned very high or very low ratings. Using a regularization coefficient $\lambda$, the estimation can be performed by solving the following unconstrained quadratic program:
\begin{align}
\min_{\beta} \left \{ \sum_{i=1}^n ( \sum_{j \in \mc{V}} x_{ij} \beta_{j} - y_{i} )^2 + \sum_{j \in \mc{V}} (\lambda \beta_j)^2 \right \} \nonumber 
\end{align}

In the context of plus-minus ratings for soccer players, let $\mc{M}$ be a set of matches. Each match $m \in \mc{M}$ can be divided in a number of segments $s \in \mc{S}_m$, where each player on the field during the segment is playing for the whole segment. One possibility is to split into segments for each substitution and for each time a player is sent off with a red card. The duration of segment $s$ of match $m$ is $d(m,s)$ minutes. For a given segment, let $f^{LHS}(m,s)$ be the left hand side of a row in the regression model, and let $f^{RHS}(m,s)$ be the right hand side. Let the regularization term for variable $j$ be denoted by $f^{REG}(\beta_j)$. Furthermore, let $w(m,s)$ be the importance of segment $s \in \mc{S}_m$ of match $m \in \mc{M}$, allowing different segments to be weighted differently when estimating the parameters of the model. Regularized plus-minus ratings can then be described by the following unconstrained quadratic program:
\begin{align}
\min Z(\beta) = \sum_{m \in \mc{M}, s \in \mc{S}_m} \left ( w(m,s) f^{LHS}(m,s) - w(m,s) f^{RHS}(m,s) \right ) ^2 + \sum_{j \in \mc{V}} \left (  f^{REG}(\beta_j) \right )^2 \nonumber
\end{align}

\noindent which by specifying the details of $f^{LHS}$, $f^{RHS}$, $f^{REG}$, and $w(m,s)$, provides a specific variant of adjusted plus-minus or regularized plus-minus ratings.

\subsection{Regularized plus-minus ratings}
\label{sec:rapm}

To obtain a plain regularized plus-minus rating, define the following. Let $h(m)$ and $a(m)$ be the two teams involved in match $m$, and let $\mc{P}_{m,s,h}$ and $\mc{P}_{m,s,a}$ be the sets of players involved in segment $s$ of match $m$ for team $h = h(m)$ and $a = a(m)$, respectively. During segment $s$ of match $m$, the number of goals scored by team $a$ and $h$ is given by $g(m,s,a)$ and $g(m,s,h)$, and the goal difference $g(m,s) = g(m,s,h) - g(m,s,a)$ is measured in favor of the home team. Then define:
\begin{align}
f^{LHS}(m,s) & = \frac{d(m,s)}{90} \left ( \sum_{p \in \mc{P}_{m,s,h}} \beta_p- \sum_{p \in \mc{P}_{m,s,a}} \beta_p \right ) \nonumber  \\
f^{RHS}(m,s) & = g(m,s) \nonumber  \\
f^{REG}(\beta_j) & = \lambda \beta_j \nonumber \\
w(m,s) & = 1 \nonumber
\end{align}


The above regularized plus-minus rating does not take into account players sent off. Hence, it seems fair to discard segments where any team does not have a full set of eleven players. This can be done by a simple redefinition of $\mc{S}_m$. A version of regularized plus-minus taking into account red cards, home advantage, and the recency of observations was presented by \citet{SaHv17}. 

\subsection{Novel regularized plus-minus ratings} \label{sec:PlayerRatings:novel}
\label{sec:pmpr}

The following describes a novel regularized plus-minus rating for football players, using an improved method to model home field advantages, an improved method to take into account red cards, and letting the ratings of players depend on their age. The method is also extended in a new way to improve the handling of players appearing in different leagues or divisions, and by introducing a more effective scheme for setting segment weights. The rating model aims to explain the observation defined as $f^{RHS}(m,s)$ using variables and parameters expressed as $f^{LHS}(m,s)$. The latter can be divided into 1) components that depend on the players $p$ involved in the segment, $f^{PLAYER}(m,s,p)$, 2) components that depend on the segment $s$ but not the players, $f^{SEGMENT}(m,s)$, and 3) components that depend on the match $m$, but not the segment, $f^{MATCH}(m)$. Thus we can write:
\begin{align}
\frac{90}{d(m,s)} f^{LHS}(m,s) = & \sum_{p \in \mc{P}_{m,s,h}} f^{PLAYER}(m,s,p) \nonumber \\
               & - \sum_{p \in \mc{P}_{m,s,a}} f^{PLAYER}(m,s,p) \nonumber \\
               & + f^{SEGMENT}(m,s) \nonumber \\
               & + f^{MATCH}(m) \nonumber
\end{align}

\noindent where setting $f^{SEGMENT}(m,s) = f^{MATCH}(m) = 0$ and $f^{PLAYER}(m,s,p) = \beta_p$ gives the plain regularized plus-minus as defined earlier. However, instead of just following the structure of a regularized regression model and making improvements to the variables included, the novel ratings also exploit that the ratings are described by an unconstrained quadratic program. In particular, the regularization terms for some of the variables are replaced by more complex expressions.


The home field advantage may vary between different league systems. For example, since the home field advantage is measured in terms of the goal difference per 90 minutes, it may be that the advantage is different in high scoring and low scoring tournaments. Let $c(m)$ be the country or competition type in which match $m$ takes place. Home field advantage is then modelled by setting
\begin{align}
  f^{MATCH}(m) & = \left \{
    \begin{array}{lr}
    \beta^H_{c(m)} & \textrm{if team~} h(m) \textrm{~has home advantage} \\
    0   & \textrm{otherwise}
    \end{array}
  \right. \nonumber
\end{align}


To correctly include the effect of players being sent off after red cards, the average rating of the players left on the pitch is used as the baseline to which additional variables corresponding to the effect of red cards on the expected goal differences are added. To this end, $f^{LHS}$ is first redefined as follows:
\begin{align}
\frac{90}{d(m,s)} f^{LHS}(m,s) = & \frac{11}{|\mc{P}_{m,s,h}|} \sum_{p \in \mc{P}_{m,s,h}} f^{PLAYER}(m,s,p) \nonumber \\
               & - \frac{11}{|\mc{P}_{m,s,a}|} \sum_{p \in \mc{P}_{m,s,a}}  f^{PLAYER}(m,s,p) \nonumber \\
               & + f^{SEGMENT}(m,s)  \nonumber\\
               & + f^{MATCH}(m)  \nonumber
\end{align}

Now, define $r(m,s,n)= 1$ if team $h$ has received $n$ red cards and team $a$ has not, $r(m,s,n)= -1$ if team $a$ has received $n$ red cards and team $h$ has not, and $r(m,s,n)= 0$ otherwise. Then, red card variables are introduced, where a difference is made between the value of a red card for the home team and for the away team, by rewriting $f^{SEGMENT}(m,s)$ as:
\begin{align}
  f^{SEGMENT}(m,s) & = \sum_{n = 1}^4 r(m,s,n) \beta^{HOMERED}_n, & \textrm{if~} \sum_{n = 1}^4 r(m,s,1) \geq 0 \nonumber \\
  f^{SEGMENT}(m,s) & = \sum_{n = 1}^4 r(m,s,n) \beta^{AWAYRED}_n, & \textrm{if~} \sum_{n = 1}^4 r(m,s,1) < 0 \nonumber
\end{align}


Playing strength is not constant throughout a player's career. In particular, being too young and inexperienced or too old and physically deteriorated, may both be seen as disadvantageous. In a paper devoted to studying the peak age of football players, \citet{De16} took performance ratings as given (calculated by a popular web page for football statistics), and fit different models to estimate the age effects. In that study, the peak age of players was estimated to between 25 and 27 years, depending on the position of the players.

To include an age effect, the player rating component of the model, $f^{PLAYER}(m,s,p)$, is modified. Let $t = t(m)$ denote the time when match $m$ is played, and let $T$ denote the time that the ratings are calculated. Let $t(m,p)$ be the age of player $p$ at time $t(m)$. The ages of players, $t(m,p)$, are measured in years. In addition to considering quadratic and cubic functions to describe the effect of a player's age, \citet{De16} introduced separate dummy variables for each age, year by year. In the regularized plus-minus model, this can be mimicked by representing the age effect as a piecewise linear function. To accomplish this, define a set of ages, $\mc{Y} = \{y^{min}, y^{min}+1, \ldots, y^{max}\}$. For a given match $m$ and player $p$, let $\max\{\min\{t(m,p),y^{max}\}, y^{min}\} = \sum_{y \in \mc{Y}} u(y,t(m),p) y$, where $\sum_{y \in \mc{Y}} u(y,t(m),p) = 1$, $0 \leq u(y,t(m),p) \leq 1$, at most two values $u(y,t(m),p)$ are non-zero, and if there are two non-zero values $u(y,t(m),p)$ they are for consecutive values of $y$. 
The player component can then be stated as:
\begin{align}
f^{PLAYER}(m,s,p) = & \beta_p + \sum_{y \in \mc{Y}} u(y,t(m),p) \beta^{AGE}_y \nonumber
\end{align}

If all players in a match have the exact same age, the age variables cancel out. However, when players are of different age, the corresponding effects of the age difference can be estimated. As a result, players are not assumed to have a constant rating over the entire time horizon of the data set, but are instead assumed to have a rating that follows an estimated age curve.

The regularization terms are not strictly necessary for variables other than the player rating variables, $\beta_p$. However, for smaller data sets, it seems beneficial to include the regularization terms also for additional variables, such as for the home field advantage and the red card effects. For the age variables, $\beta^{AGE}_y$, a different scheme is chosen, as it seems beneficial to make sure that the estimates for each age are somehow smoothed. This can be accomplished by the following replacements for the regularization terms:
\begin{align}
f^{REG}(\beta^{AGE}_y) & = \lambda \left(\beta^{AGE}_y - (\beta^{AGE}_{y-1} + \beta^{AGE}_{y+1})/2 \right), & y \in \mc{Y} \setminus \{y^{min}, y^{max}\} \nonumber \\
f^{REG}(\beta^{AGE}_y) & = 0, & y \in \{y^{min}, y^{max}\} \nonumber
\end{align}


For players with few minutes of recorded playing time, the standard regularization ensures that the players' ratings are close to zero. \citet{SaHv17} included a tournament factor in the player ratings, thus allowing players making their debut in a high level league to obtain a higher rating than players making their debut in low level leagues. This tournament factor is generalized here, as follows. Let $\mc{B}$ be a set of different leagues, and let $\mc{B}_p \subseteq \mc{B}$ be the set of leagues in which player $p$ has participated. The player component is then further refined to become
\begin{align}
f^{PLAYER}(m,s,p) = & \beta_p + \sum_{y \in \mc{Y}} u(y,t(m),p) \beta^{AGE}_y + \frac{1}{|\mc{B}_p|} \sum_{b \in \mc{B}_p} \beta^{B}_b \nonumber
\end{align}

This helps to discriminate players from different leagues. However, a further refinement of this is achieved by modifying the regularization terms. Instead of always shrinking a player's individual rating component $\beta_p$ towards 0, as in the plain regularized plus-minus ratings, the whole expression providing the current rating of a player is shrunk towards a value that depends on a set of similar players. Let $\mc{P}^{SIMILAR}_p$ be a set of players that are assumed to be similar to player $p$. In this work, the set is established by using the teammates of $p$ that have been on the pitch together with $p$ for the highest number of minutes. Let $t(p, p \prime)$ be the time of the last match where players $p$ and $p \prime$ appeared on the pitch for the same team. Now, define the following auxiliary expression, where $w^{AGE}$ is a weight for the influence of the age factor:
\begin{align}
f^{AUX}(p,t,w^{AGE}) = & \beta_p + w^{AGE} \sum_{y \in \mc{Y}} u(y,t,p) \beta^{AGE}_y + \frac{1}{|\mc{B}_p|} \sum_{b \in \mc{B}_p} \beta^{B}_b \nonumber
\end{align}

The rating of player $p$ at time $T$ is then equal to $f^{AUX}(p,T,1)$, and it is this value that will be shrunk towards a value that depends on the teammates of $p$, rather than towards 0. To this end, the regularization term for player $p$ is replaced by the following:
\begin{align}
f^{REG}(\beta_p) & = \lambda \left ( f^{AUX}(p,T,1) - \frac{ w^{SIMILAR} }{|\mc{P}^{SIMILAR}_p|} \left ( \sum_{p \prime \in \mc{P}^{SIMILAR}_p} f^{AUX} (p \prime, t(p,p \prime), w^{AGE}) \right ) \right ) \nonumber
\end{align}

\noindent where $w^{SIMILAR} \leq 1$ is another weight that controls the emphasis of shrinking the rating of player $p$ towards the rating of similar players versus shrinking towards 0. This replacement of the regularization terms for player rating components, together with the modified regularization terms for the age components, makes the full model incompatible with the framework of regularized linear regression models. Instead, the model is interpreted as an unconstrained quadratic program.


The model estimation is performed by minimizing the sum of squared deviations between observed goal differences and a linear expression of player ratings and additional factors. The sum is taken over all segments from all matches included in the data. However, not all of these segments are equally informative, and better ratings can be obtained by changing the relative weight $w(m,s)$ of different segments. 

The weights used here have three components. The first component emphasizes that more recent matches are more representative for the current strength of players. Hence, a factor $w^{TIME}(m) = e^{\rho_1(T-t(m))}$ is included, which leads to smaller weights for older matches. The second component focuses on the duration of a segment, with longer segments being more important than shorter segments. Given two parameters $\rho_2$ and $\rho_3$, and the duration of a segment, $d(m,s)$, a factor on the form $w^{DURATION}(m,s) = (d(m,s) + \rho_2) / \rho_3$ is included. The third component takes into account the goal difference at the beginning of the segment, $g^0(m,s)$, as well as the goal difference at the end of the segment, $g^1(m,s) = g^0(m,s)+g(m,s)$, introducing the factor $w^{GOALS}(m,s) = \rho_{4}$ if $|g^0(m,s)| \geq 2$ and $|g^1(m,s)| \geq 2$, and $w^{GOALS}(m,s) = 1$ otherwise. The weight of a segment can then be stated as:
\begin{align}
w(m,s) = w^{TIME}(m,s) w^{DURATION}(m,s) w^{GOALS}(m,s). \nonumber
\end{align}

\section{Case Studies}
\label{sec:casestudies}
In this section we describe a number of case studies used to test model \eqref{eq:FTCP}. The case studies consist of the 20 clubs competing in the English Premier League (EPL) during the 2013/14 season. Each club is characterized by the current team composition and a list of target players, and we use model \eqref{eq:FTCP} to address the transfer market of summer 2014, in preparation for season 2014/15. The data of the case studies is made available online at \url{https://github.com/GioPan/instancesFTCP}.

In \Cref{subsec:Instances} we describe the clubs and their current and target players. In \Cref{subsec:regression} we  introduce a model of the market value of the player which allows us to obtain an empirical probability distribution. Given the complexity of solving model \eqref{eq:FTCP} with the original empirical distribution, in \Cref{subsec:model:saa} we introduce its Sample Average Approximation. In \Cref{subsec:ratings} we provide some statistics about the ratings of the players in the case studies. The case studies are subsequently used to perform a number of tests which will be thoroughly described in \Cref{sec:results}.

\subsection{Clubs and players}
\label{subsec:Instances}

The case studies used for testing are adapted from those introduced by \citet{Pan17} based on the English Premier League (EPL). The case studies describe the 20 teams in the EPL 2013/2014 dealing with the summer 2014 transfer market. Each team is characterized by a budget, a set of players currently owned and the set of target players. Given a focal team among the 20 available, the set $\mc{P}$ consists of the set of current player and the set of target players. Each player is characterized by age, role, current value, purchase and sale price, loan fees, and whether the player can be purchased, sold, or temporarily change club in a loan agreement. 

In addition to the above mentioned data, we set $N=25$ in accordance with EPL rules. Furthermore, we test different formations, where a formation determines the number of players required for each role. Thus, for each role $r\in \mc{R}$ we set $\underline{N}_r$ according to \Cref{tab:cs:formations}, and $\overline{N}_r=\infty$, implying that it is allowed to have more than $\underline{N}_r$ players covering role $r$. The roles used here are simply player positions: goalkeeper (GK), right-back (RB), centre-back (CB), left-back (LB), right winger (RW), centre midfielder (CM), left winger (LW), attacking midfielder (AM), and forward (FW). Finally, we set $V$ equal to the initial market value of the team (i.e., the club wishes to ensure a non-decreasing value of the team) and we use a $7\%$ discount factor. Uncertain values and ratings are discussed in \Cref{subsec:regression} and \Cref{subsec:ratings}, respectively. 

\begin{table}[htb]
  \centering
  \caption{Formations and players required in each role ($\underbar{N}_r$). }
  \label{tab:cs:formations}
  \begin{tabular}{cccccccccc}
    \toprule
              & \multicolumn{9}{|c}{$r\in\mc{R}$}      \\
    Formation & GK & RB & CB & LB & RW & CM & LW & AM & FW \\
    \midrule
    442       & 3  & 2  & 4  & 2  & 2  & 4  & 2  &0 &4  \\
    433       & 3  & 2  & 4  & 2  & 0  & 6  & 0  &0 &6  \\
    4312       & 3  & 2  & 4  & 2  & 0  & 6  & 0  &2 &4  \\
    352       & 3  & 0  & 6  & 0  & 2  & 6  & 2  &0 &4  \\
    343       & 3  & 0  & 6  & 0  & 2  & 4  & 2  &0 &6  \\
    \bottomrule
  \end{tabular}
\end{table}

\subsection{Modeling the Uncertainty}
\label{subsec:regression}

Uncertain future player values are modeled using a linear regression model whose parameters were estimated using English Premier League data for the seasons 2011/2012 through 2016/2017. Consistently with \citet{Pan17} we use the current value of the player and its role as
explanatory variables and the value after one season as the dependent variable. Unlike in \citet{Pan17}, different
regression models are used for different age intervals, as this allows us to capture the higher volatility in the value of younger players. 
In this exercise, we aimed to obtain a sufficiently good model of the available data, indicated by high values of the $R^2$ value and no evidence of non-linear patterns in the residual analysis. Therefore, the resulting model should not be understood as a good prediction tool, as no out-of-sample analysis was performed.

For each player $p\in \mc{P}$, and scenario $s\in\mc{S}$, the future value
$V_{ps}$ is thus obtained as
\begin{equation}
  \label{eq:cs:regressionmodel}
  V_{ps}=\bigg(\alpha_{a}\sqrt[4]{V_{p}^C}+\sum_{r\in\mc{R}}\beta_{ar}\delta(p,\mc{P}_r)\bigg)^4*(1+\epsilon_{as})
\end{equation}
where  $V^C_p$ is the current value of player $p$, $\delta(p,\mc{P}_{r})$
is an indicator function of the membership of player $p$ in $\mc{P}_r$, i.e., it is equal to $1$ if player
$p$ has role $r$, $0$ otherwise, and $\epsilon_{as}$ is an i.i.d. sample from the empirical prediction error distribution of the regression model for the specific age group. Notice that $\alpha_a$, $\beta_{ar}$ and $\epsilon_a$ are estimated for different age
intervals $a\in \{(\cdot,20],[21,22],[23,24],[25,26],[27,28],[29,30],[31,32],[33,\cdot)\}$.  While model \eqref{eq:cs:regressionmodel} does not include a constant term, the indicator function is used to determine a role-specific constant term. \Cref{tab:R2} reports the $R^2$ coefficient for the regression model in each age group.

\begin{table}[h!]

    \centering
    \caption{$R^2$ for regression model \eqref{eq:cs:regressionmodel} in each age group.}
    \label{tab:R2}
    \begin{tabular}{cc}
    \toprule
    Age group&$R^2$\\
    \midrule
    $(\cdot,20]$&$0.97$\\
    $[21,22]$&$0.98$\\
    $[23,24]$&$0.98$\\
    $[25,26]$&$0.99$\\
    $[27,28]$&$0.99$\\
    $[29,30]$&$0.99$\\
    $[31,32]$&$0.98$\\
    $[33,\cdot)$&$0.98$  \\
    \bottomrule
    \end{tabular}
\end{table}

Not all roles were statistically significant as predictors ($p$-value smaller than $0.05$). Particularly, the roles which were statistically insignificant for some age group are reported in \Cref{tab:pvalues}.

\begin{table}[h!]

    \centering
    \caption{Roles statistically insignificant as predictors.}
    \label{tab:pvalues}
    \begin{tabular}{cc}
    \toprule
    Role     & Age group in which the role is statistically insignificant as a predictor \\
    \midrule
        Secondary striker & $(0,20]$, $[29,30]$ \\
        Left midfielder & $[21,22]$, $[27,28]$ ,$[29,30]$ $[31,32]$\\
        Right midfielder & $[23,24]$, $[25,26]$, $[27,28]$, $[31,32]$\\ 
        Left wing &$[31,32]$\\
         \bottomrule
    \end{tabular}
\end{table}

\subsection{Sample Average Approximation}
\label{subsec:model:saa}

The empirical prediction error distribution for regression model \eqref{eq:cs:regressionmodel} is a discrete distribution with a very large support. This makes model \eqref{eq:FTCP} a mixed-integer linear program. However, solving the model using directly the original discrete distribution is impractical due to the large number of realizations of the player values. Therefore, we solve its \textit{Sample Average Approximation} (SAA), see e.g., \citet{KleS02,Sha03,PagAS09}. 

Let $\mc{S}=\{1,\ldots,S\}$ and let $(V_{ps})_{s\in \mc{S}}$ be an $|\mc{S}|$-dimensional i.i.d. sample of $\tilde{V}_p$ for all $p\in \mc{P}$. Furthermore, for all $s\in \mc{S}$ let $w_s$ be a binary variable which is equal to $1$ if the team value exceeds the threshold $V$ under scenario $s$, $0$ otherwise. That is $w_s = 1 \implies \sum_{p\in \mc{P}}V_{ps}y_p \geq V$.
Constraint \eqref{eq:chance} can be approximated by constraints
\eqref{eq:chanceSAA1}-\eqref{eq:chanceSAA3}
\begin{align}
  \label{eq:chanceSAA1}  &\sum_{p\in \mc{P}}V_{ps}y_p + M(1-w_s) \geq V &~~ s\in \mc{S}\\
  \label{eq:chanceSAA2}  &\sum_{s\in \mc{S}}\frac{1}{|\mc{S}|}w_s \geq \alpha &\\
  \label{eq:chanceSAA3}  & w_s \in \{0,1\} & ~~s\in \mc{S}
\end{align}
where $M$ is a suitable upper bound to $V-\sum_{p\in \mc{P}}V_{ps}y_p$. The quantity $V-\sum_{p\in \mc{P}}V_{ps}y_p$ is bounded above by $M=V$, achieved when $\sum_{p\in \mc{P}}V_{ps}y_p = 0$ (which yields an infeasible solution).

We approximate the original empirical distribution by means of a sample of size $70$. Numerical tests show that this sample size ensures both in-sample stability of the objective function (i.e., a negligible standard deviation of the SAA optimal objective value across different samples of size $70$) and out-of-sample satisfaction of the chance constraint, assessed on a sample of $1000$ scenarios.
That is, we observe that the change of optimal objective value determined by different samples of size $70$ is relatively small. In \Cref{tab:iss}, which reports the results for the cases with $\alpha=0.2$ and $0.8$, we notice that the standard deviation of the optimal objective value is typically at least two orders of magnitude smaller than the mean. We also observe that the solutions which satisfy the chance constraint with an i.i.d. sample of size $70$, also satisfy the chance constraint with a sample of size $1000$. \Cref{tab:oos} reports some statistics about the out-of-sample probability obtained with $\alpha=0.2$ and $0.8$. It can be seen that the out-of-sample probability is always greater than $\alpha$. Stability tests have been run with formation 443 which, as explained later, is used throughout the computational study.

\begin{table}[htb]

    \centering
        \caption{Results of the in-sample stability test for SAA obtained with samples of size $70$ and $\alpha=0.2$ and $0.8$. The mean and standard deviation of the objective values were calculated by solving $10$ different SAAs, each obtained with a sample of size $70$.}
    \label{tab:iss}
    \begin{tabular}{c|cccc}
    \toprule
    &   \multicolumn{2}{c}{$\alpha=0.2$} & \multicolumn{2}{c}{$\alpha=0.8$}\\
        Team    & Avg. obj. value & St. dev & Avg. obj. value & St. dev\\
        \midrule
Arsenal FC &                 $4.02$ & $ 0.00 $ &    $ 4.01 $ & $7.72 \times 10^{-3}$\\
Aston Villa &               $2.68$  & $1.98 \times 10^{-4} $   &   $2.67$ & $ 3.81 \times 10^{-3}$\\
Cardiff City &              $2.08 $ & $0.00 $      & $2.05 $ & $6.11 \times 10^{-4}$\\
Chelsea FC    &             $4.21$  &$ 8.56 \times 10^{-3} $&    $  4.03$ &  $2.74 \times 10^{-2}$\\
Crystal Palace &            $1.95 $ & $2.34 \times 10^{-16}$ &    $  1.95 $&  $2.34 \times 10^{-16}$\\
Everton FC      &           $2.81 $ & $4.68 \times 10^{-16}$ &      $2.81$ &  $4.68 \times 10^{-16}$\\
Fulham FC        &          $2.45 $ & $1.73 \times 10^{-3}$ &      $2.45$ & $ 2.82 \times 10^{-3}$\\
Hull City        &          $2.18$ & $0.00  $     &  $2.18 $ & $0.00$\\
Liverpool FC     &          $3.74 $ & $1.49 \times 10^{-3}$  &     $3.74 $ &$ 0.00$\\
Manchester City  &          $4.50 $ &   $6.59 \times 10^{-3} $  &    $4.016$  & $1.28 \times 10^{-1}$\\
Manchester United&          $4.40 $ & $3.54 \times 10^{-3} $&      $ 4.40$ & $ 6.62 \times 10^{-16}$\\
Newcastle United  &         $2.22$  & $4.31 \times 10^{-3}$  &    $ 2.22$ & $ 0.00$\\
Norwich City      &         $2.37$ & $4.22 \times 10^{-3} $  &    $2.22$ &  $2.69 \times 10^{-2}$\\
Southampton FC    &         $2.74$ & $0.00  $    & $2.69 $ & $6.86 \times 10^{-3}$\\
Stoke City        &         $2.16 $ & $ 4.68 \times 10^{-16}$    &  $ 2.14$& $  5.79 \times 10^{-3}$\\
Sunderland AFC    &         $2.26 $ & $0.00  $    &$ 2.26$   &$0.00$\\
Swansea City      &         $2.70 $ & $ 8.28 \times 10^{-4}$    &   $2.67$&   $8.14 \times 10^{-3}$\\
Tottenham Hotspur  &        $3.03 $ & $0.00$    &   $2.90$  & $2.09 \times 10^{-2}$\\
West Bromwich Albion&       $2.18 $ & $0.00$    &  $ 2.18 $ & $0.00$\\
West Ham United      &      $2.23$  & $4.68 \times 10^{-16}$    &  $2.23$   &$4.68 \times 10^{-16}$\\
\bottomrule
    \end{tabular}
\end{table}

\begin{table}[htb]

    \centering
    \caption{Statistics on the out-of-sample probability of satisfying the $\alpha=0.2$ and $\alpha=0.8$ chance constraint, measured with a sample of size $1000$. The tests were run with formation $433$.}
    \label{tab:oos}
    \begin{tabular}{c|cccc}
    \toprule
    & \multicolumn{4}{c}{Out-of-sample probability}\\
        $\alpha$ & Min   & Average & St. dev. & Max  \\
        \midrule
        $0.2$ &0.27 & 0.79 & 0.23 & 1.00\\
        $0.8$ & 0.80 & 0.94& 0.04 &1.00\\ 
        \bottomrule
    \end{tabular}
\end{table}

\begin{figure}[h!]
    \hspace*{-4cm}
    \includegraphics[width=1.4\textwidth]{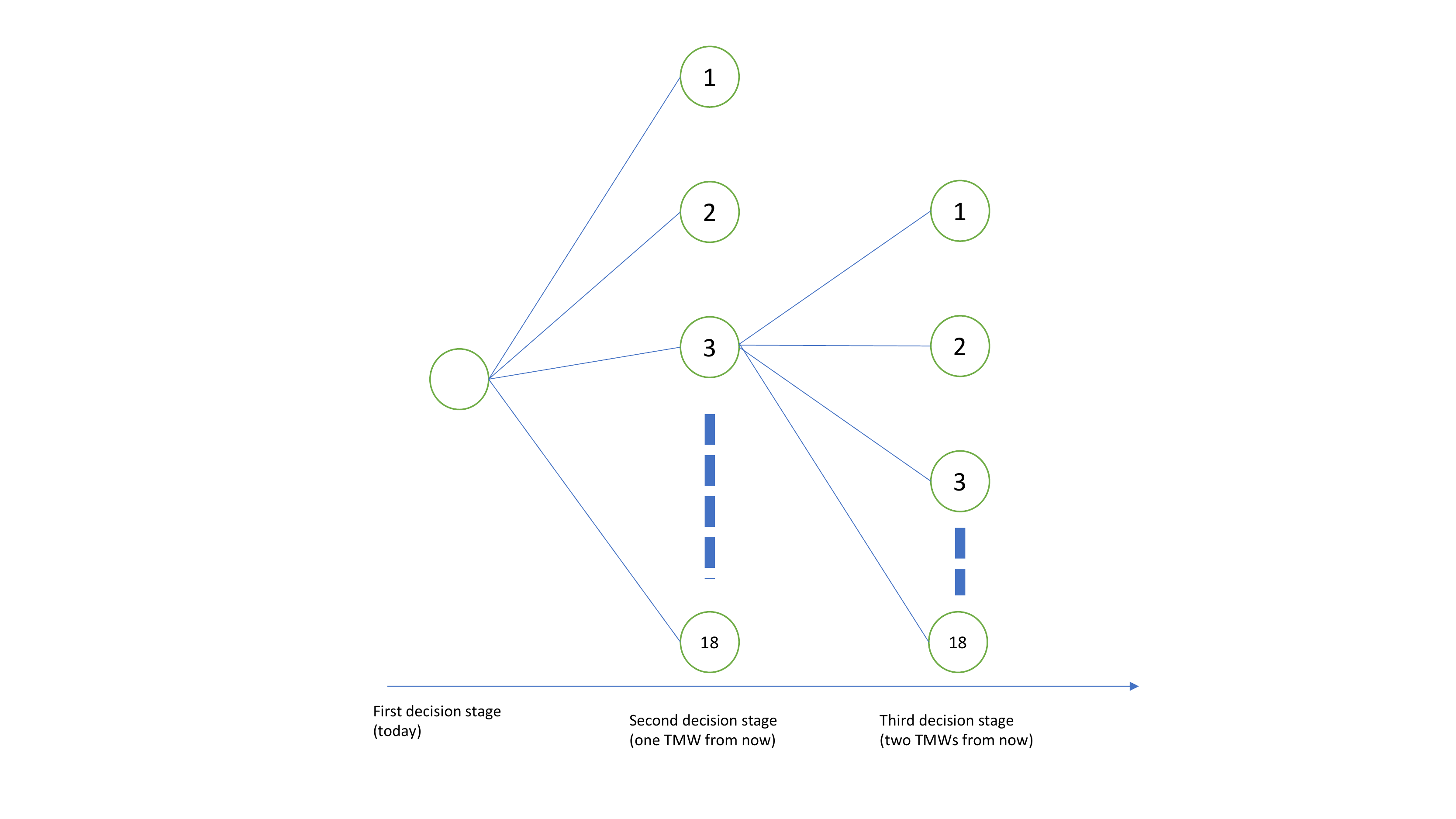}
    \caption{Qualitative description of a three-stage scenario tree. Each node represents a possible realization of the player value at the corresponding TMW, and is obtained using the age, role, and value at the previous TMW, according to model \eqref{eq:cs:regressionmodel}.}
    \label{fig:tree}
\end{figure}

When solving the multistage stochastic program introduced by \citet{Pan17}, scenario trees are also sampled from the empirical prediction error distribution based on regression model \eqref{eq:cs:regressionmodel}. The in-sample stability test according to \cite{KauW07} was performed by \citet{Pan17}, who report that the model was in-sample stable when drawing 18 conditional realizations per stage. 
That is, given current player information (age, value, and role), the market value distribution at the next TMW (second stage) is approximated by 18 realizations. Following, for each of the 18 realizations in the second stage we draw 18 conditional realizations for the third stage, and so on. The procedure is sketched in \Cref{fig:tree}. We adopt the same scenario tree structure.

\clearpage

\subsection{Ratings}
\label{subsec:ratings}

The player ratings used in this study are calculated using the model outlined in \Cref{sec:PlayerRatings:novel}. Data from more than 84,000 matches to be used in the calculations were collected from online sources. The matches come from national leagues of 25 different countries, as well as from international tournaments for club teams and national teams. There are in total 60,484 players in the data set. When calculating ratings as of July 1 2014, the unconstrained quadratic program has 60,592 variables and 598,697 squared terms. \Cref{tab:top10} shows the ten highest ranked players as of July 1 2014, consider only players with at least one match played during the last year.

\begin{table}[h!]
  \centering
  \caption{Ten highest ranked players at July 1 2014.}
    \begin{tabular}{rrccrrr}
    \toprule
    \multicolumn{1}{c}{Rank} & \multicolumn{1}{c}{Player} & Age   & Position & \multicolumn{1}{c}{Nationality} & \multicolumn{1}{c}{Team} & \multicolumn{1}{c}{Rating} \\
    \midrule
    1     & Cristiano Ronaldo & 29    & FW     & Portugal & Real Madrid & 0.317 \\
    2     & Ga\"{e}l Clichy & 28    & LB     & France & Manchester City & 0.296 \\
    3     & Lionel Messi & 27    & FW    & Argentina & Barcelona & 0.286 \\
    4     & Karim Benzema & 26    & FW     & France & Real Madrid & 0.286 \\
    5     & Thomas M\"{u}ller & 24    & FW     & Germany & Bayern M\"{u}nchen & 0.281 \\
    6     & Mesut \"{O}zil & 25    & AM    & Germany & Arsenal & 0.277 \\
    7     & Arjen Robben & 30    & FW     & Netherlands & Bayern M\"{u}nchen & 0.270 \\
    8     & J\'{e}r\^{o}me Boateng & 25    & LB    & Germany & Bayern M\"{u}nchen & 0.270 \\
    9     & Cesc F\`{a}bregas & 27    & CM    & Spain & Barcelona & 0.269 \\
    10    & Marcelo & 26    & RB    & Brazil & Real Madrid & 0.267 \\
    \bottomrule
    \end{tabular}%
  \label{tab:top10}%
\end{table}%

Parameter values for the rating calculations were determined using a different data set, containing more recent results but much fewer leagues. Using an ordered logit regression model to predict match results based on the difference in average player ratings for the two teams involved, parameters were set to minimize the quadratic loss of predictions from out-of-sample matches. This resulted in the following parameter values. The age variables are defined for $\mc{Y} = \{y^{min}=16, \ldots, y^{max}=42\}$. Observations are discounted over time with a factor of $\rho_1 = 0.1$, and are weighted using $\rho_2 = 300.0$, $\rho_3 = 300.0$, and $\rho_{4} = 2.5$. The general regularization parameter is $\lambda = 16.0$, and to make sure each player is assumed to be somewhat similar to the most common teammates, we end up with $w^{SIMILAR} = 0.85$, $w^{AGE} = 0.35$, and with the maximum number of teammates considered being 35 (for $|\mc{P}^{SIMILAR}_p|$). Deviating from these parameter settings led to worse predictions for match outcomes on the selected out-of-sample matches. 

\begin{figure}[htb]
\centering
\includegraphics[width=0.9\textwidth]{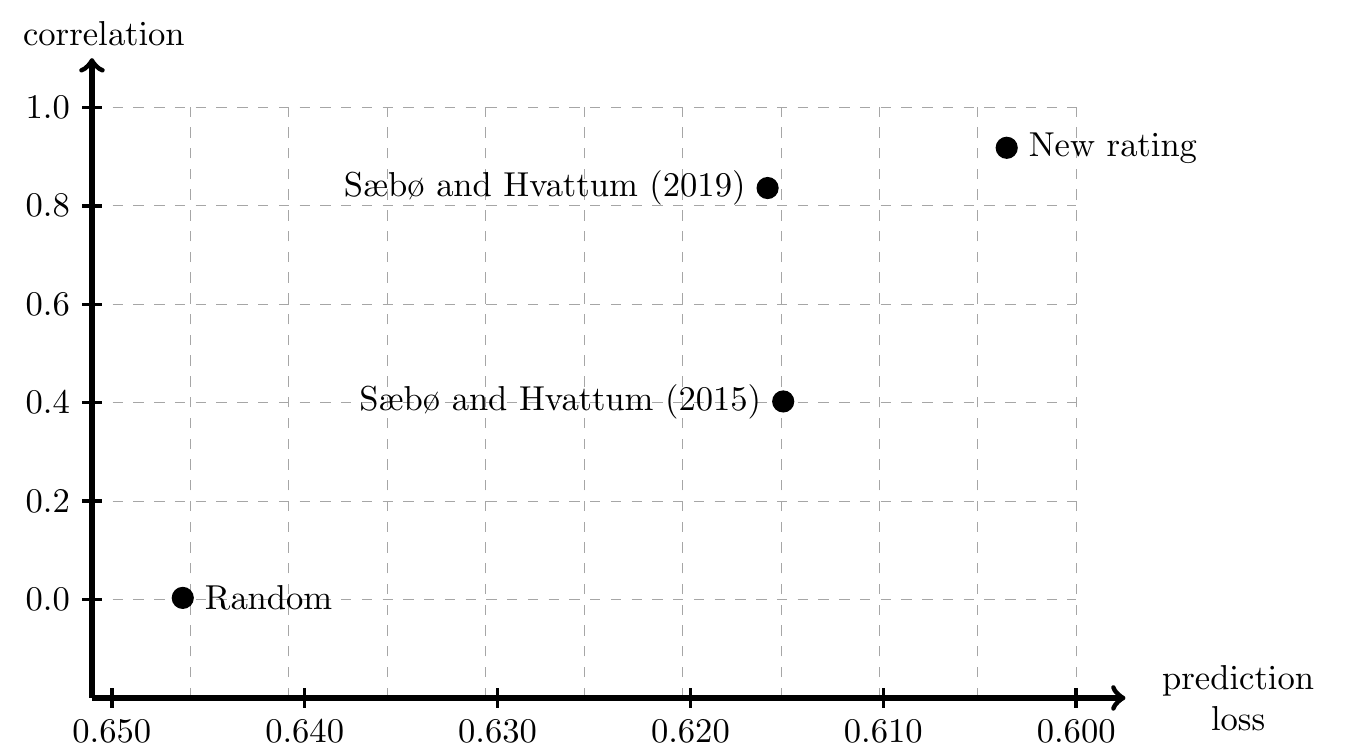}
\caption{Evaluation of the new plus-minus ratings compared to previously published ratings. The axes are oriented so that better values are found in the upper right.}\label{fig:evaluation}
\end{figure}

The new rating system is compared to two previous versions of regularized adjusted plus-minus ratings and a naive benchmark in \Cref{fig:evaluation}. The evaluation is performed along two axes: the first is the average quadratic loss on 13,800 match forecasts based on an ordered logit regression model, and the second is the Pearson correlation coefficient for ratings calculated after randomly splitting the training data in two halves, averaged over twenty repetitions. While the former represents the validity of the ratings, with lower prediction loss indicating more meaningful ratings, the latter represents the reliability of the ratings, with higher correlation indicating a more accurate rating.

\begin{figure}[htb]
\centering
\includegraphics[width=0.8\textwidth]{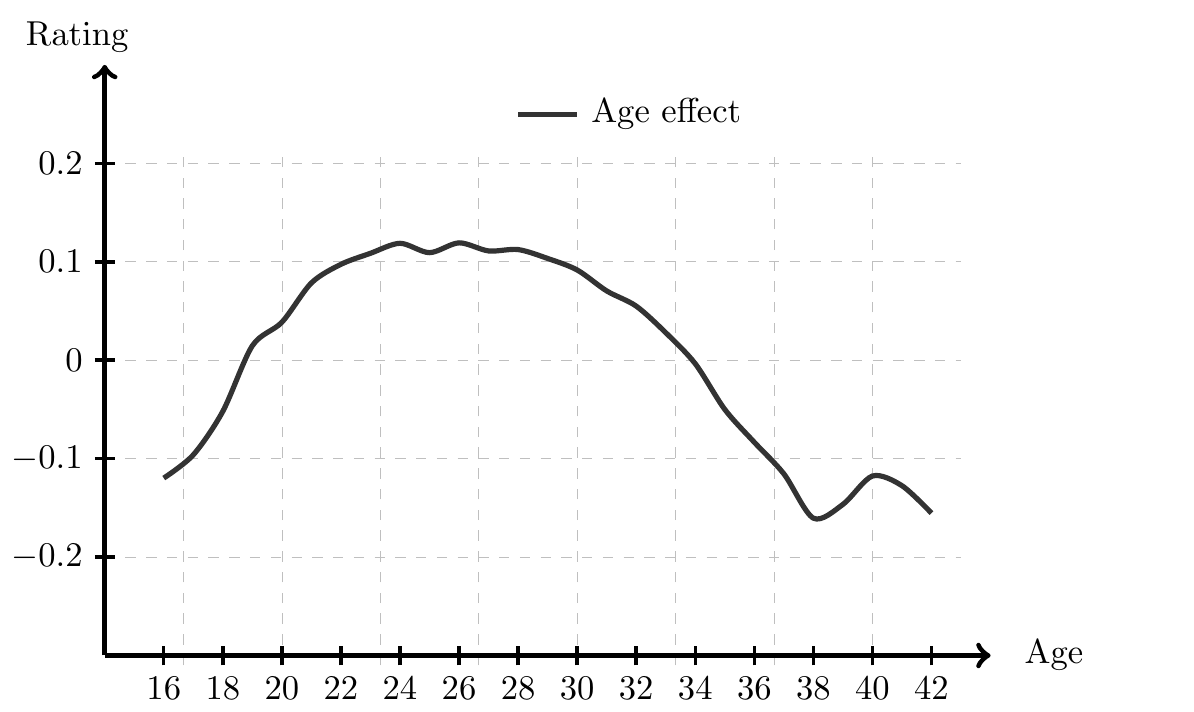}
\caption{The effect of age on player ratings, as estimated by the rating model.}\label{fig:ageProfile}
\end{figure}

\Cref{fig:ageProfile} shows how the rating model is estimating the effect of a player's age on his performance. As in \citep{De16}, it is found that the peak age is around 25--27 years. There are few observations with players aged above 40 years, and in combination with a survival bias, the estimated age curve is unreliable for relatively old players.

The rating model also provides estimates for the value of the home field advantage and the effect of red cards. The home field advantage is allowed to vary between countries, and the average effect corresponds to 0.25 goals per 90 minutes. The home field advantage is highest in the Champions League and the Europa League, with 0.41 goals per 90 minutes. The values of red cards differ between home teams and away teams in the rating model. The first red card is worth more when the away team has a player sent off. In that case, the advantage for the home team is 1.07 goals per 90 minutes, whereas the effect is 0.83 goals per 90 minutes when the home team is reduced by one man. As in \citep{SaHv17}, it is found that subsequent red cards have smaller consequences.

\section{Results and discussion}
\label{sec:results}
In this section we present the results of a number of tests performed on the case studies based on the EPL 2013/14. The scope of the tests is to illustrate the team composition strategies obtainable with model \eqref{eq:FTCP}, and particularly compare those with the strategies of a club that maximizes the team value. A team value-maximizer is modeled by means of the multistage stochastic program from \citet{Pan17}.
Note that it is expected that the two models provide different team composition strategies. Therefore, the scope of the computational study is that of examining such differences.
Furthermore, we assess the impact of different financial risk tolerances of the clubs. Unless otherwise specified we show the results obtained using a 4-4-3 formation. In \Cref{app:formation} we show that our findings are to a large extent insensitive to the formation chosen. In what follows, we refer to the chance constrained model \eqref{eq:FTCP} as CC and to the multistage stochastic program from \citet{Pan17} as MSP.

\subsection{Maximizing team value vs maximizing ratings}\label{sec:results:risk}
The solutions to the CC model are compared to those obtained by solving the MSP model. This corresponds to comparing the maximization of the ratings (subject to probabilistic constraints on the market value) to the maximization of market values regardless of player ratings. For the MSP model we consider three stages and generate 18 conditional realizations at each stage as in \citep{Pan17}, resulting in 324 scenarios (see \Cref{subsec:model:saa} for details on the scenario generation).

\begin{figure}[htb]
\includegraphics[width=\textwidth]{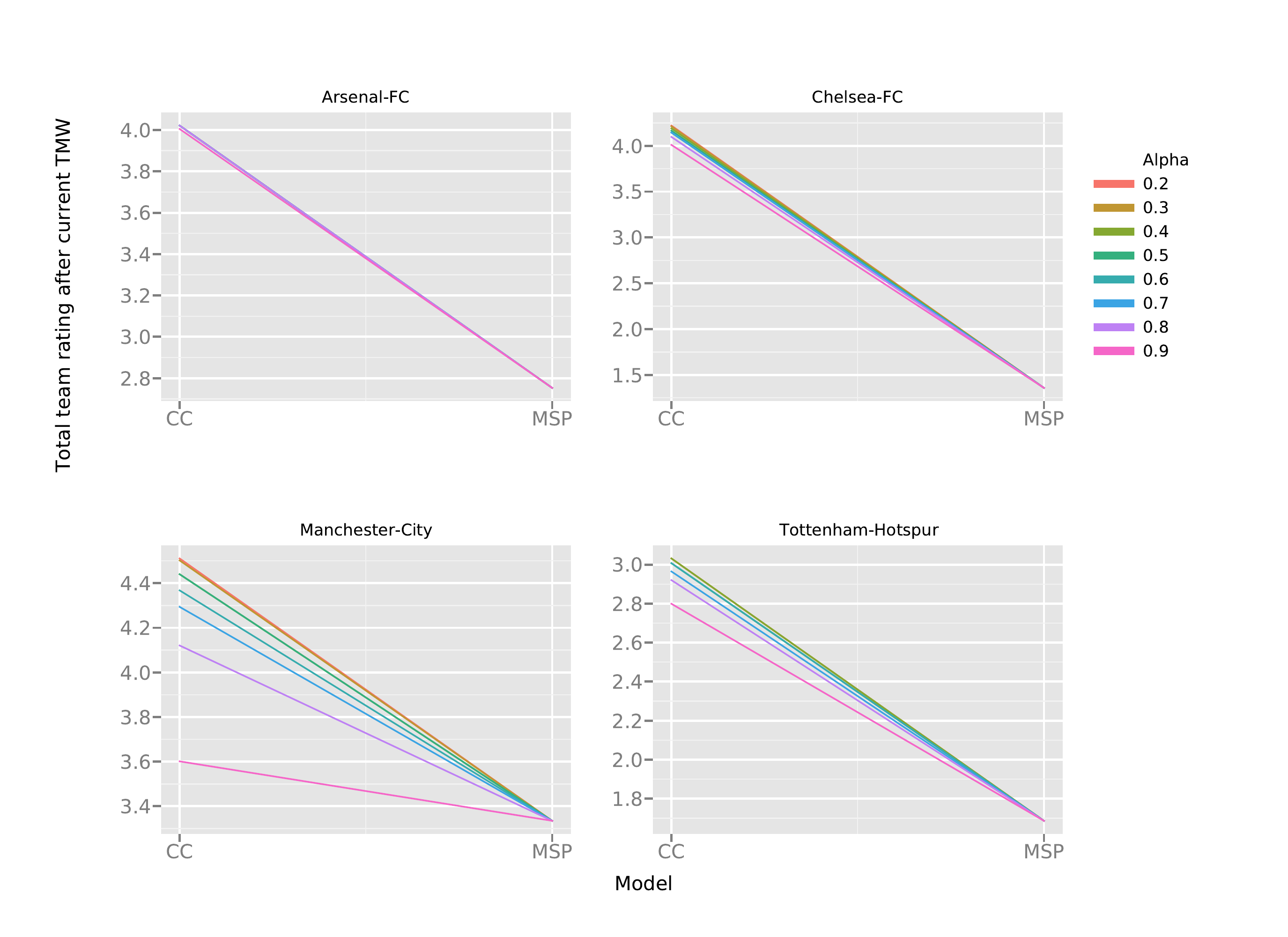}
\caption{Total rating of the team composed by the MSP and by the CC model with different levels of $\alpha$.}\label{fig:ccVSmspRating}
\end{figure}

\Cref{fig:ccVSmspRating} reports, for a sample of four clubs, the total rating of the teams composed by the MSP model and by the CC model with different values of $\alpha$. The same findings apply for the teams not shown in the figure. The rating for the MSP model is calculated for the team provided after the first TMW. The rating of the team obtained by the CC model is consistently higher than the rating of the MSP, and in most cases significantly higher. The MSP model does not find any value in signing top performers per se. Rather, the MSP model looks for players whose value is likely to increase in the future as a consequence of their age, role, and current evaluation. Very often, these players are not yet top performers. 
On the other hand, the CC model looks primarily for top rated players, that is players whose performances have provided a solid contribution to their respective team's victories in past matches. 
For several clubs, similar to the case of Arsenal-FC in \Cref{fig:ccVSmspRating}, the rating is insensitive to the value of $\alpha$. This issue is properly discussed in \Cref{sec:results:ratingvstolerance}.

\begin{figure}[htb]
\includegraphics[width=\textwidth]{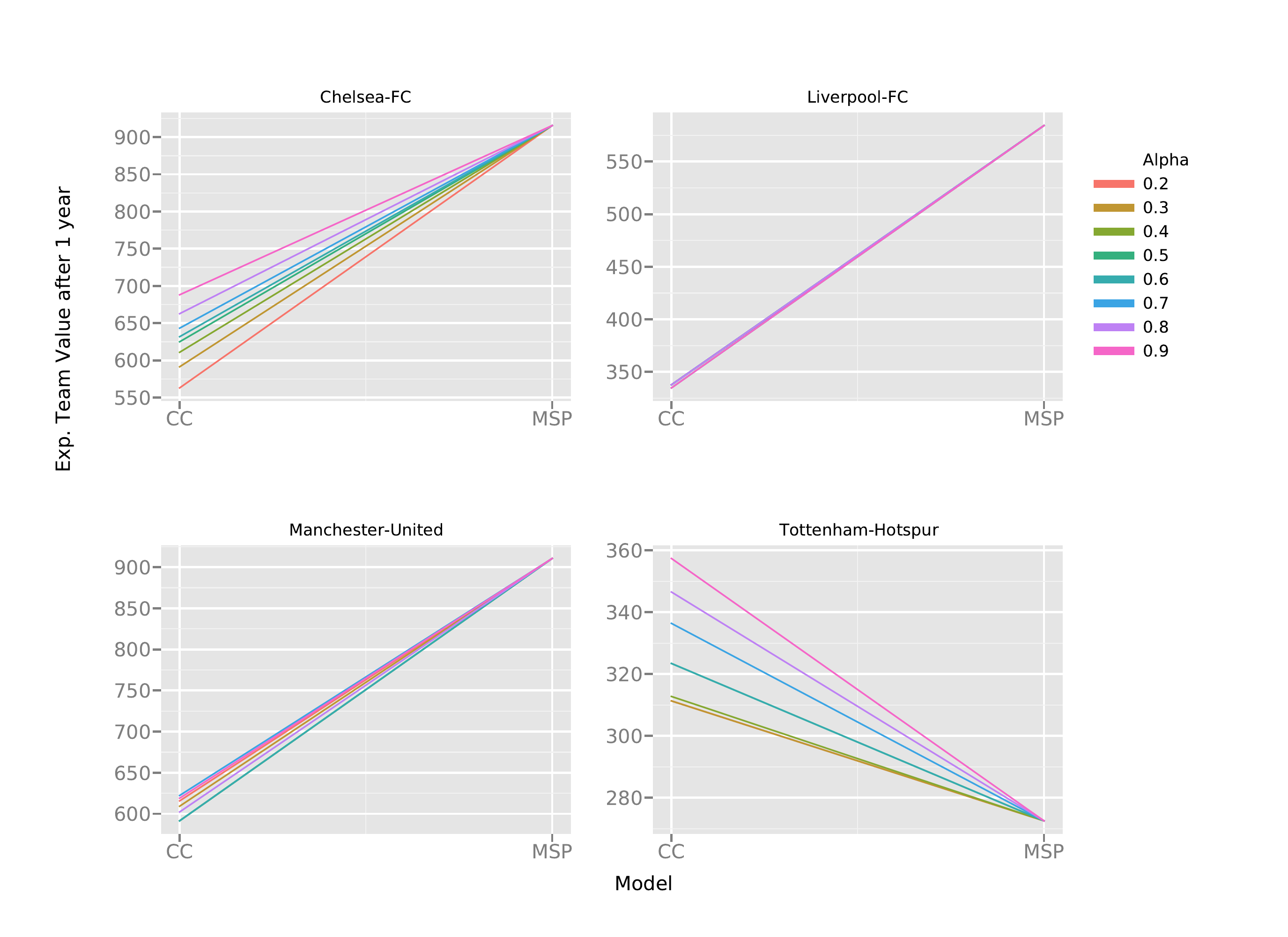}
\caption{Expected value after one year of the team composed maximizing rating and the team composed maximizing profits.}\label{fig:ccVSmspValue}
\end{figure}

Let us now turn our attention to the expected market value of the teams provided by the two models after one season, reported, for a sample of four clubs, in \Cref{fig:ccVSmspValue}. We can observe that the MSP model yields a significantly higher expected team value after one season. This is to be expected since the MSP maximizes market values. On the other hand, the CC model simply ensures that the market value of the team does not decrease after one season. Thus, the decision maker does not seek a return on the capital employed in the team, but simply wants to ensure that the investment keeps its value. In \Cref{fig:ccVSmspValue} we can also find a case for which the CC model provides a higher one-year expected team value than the MSP model (see Tottenham-Hotspur). This is due to the fact that the MSP model maximizes the average expected team value over a three-year period. Therefore, it is possible that the model suggests investments that do not necessarily yield the highest team value after one season, as long as the average over three seasons is maximized. 

\subsection{Team rating and risk tolerance}\label{sec:results:ratingvstolerance}
We illustrate the impact of the risk tolerance $\alpha$ on team ratings. We consider both the standard case in which the club wants to ensure a non-decreasing team value, and the case in which the club wants to ensure a growth of the value of the team of either $10$, $20$, or $30\%$. This corresponds to multiplying the constant $V$ by a factor $R=1.1$, $1.2$, and $1.3$, respectively, in constraint \eqref{eq:chance}. In the default case, $V$ represents the initial market value of the team.

\begin{figure}[htb]
\includegraphics[width=\textwidth]{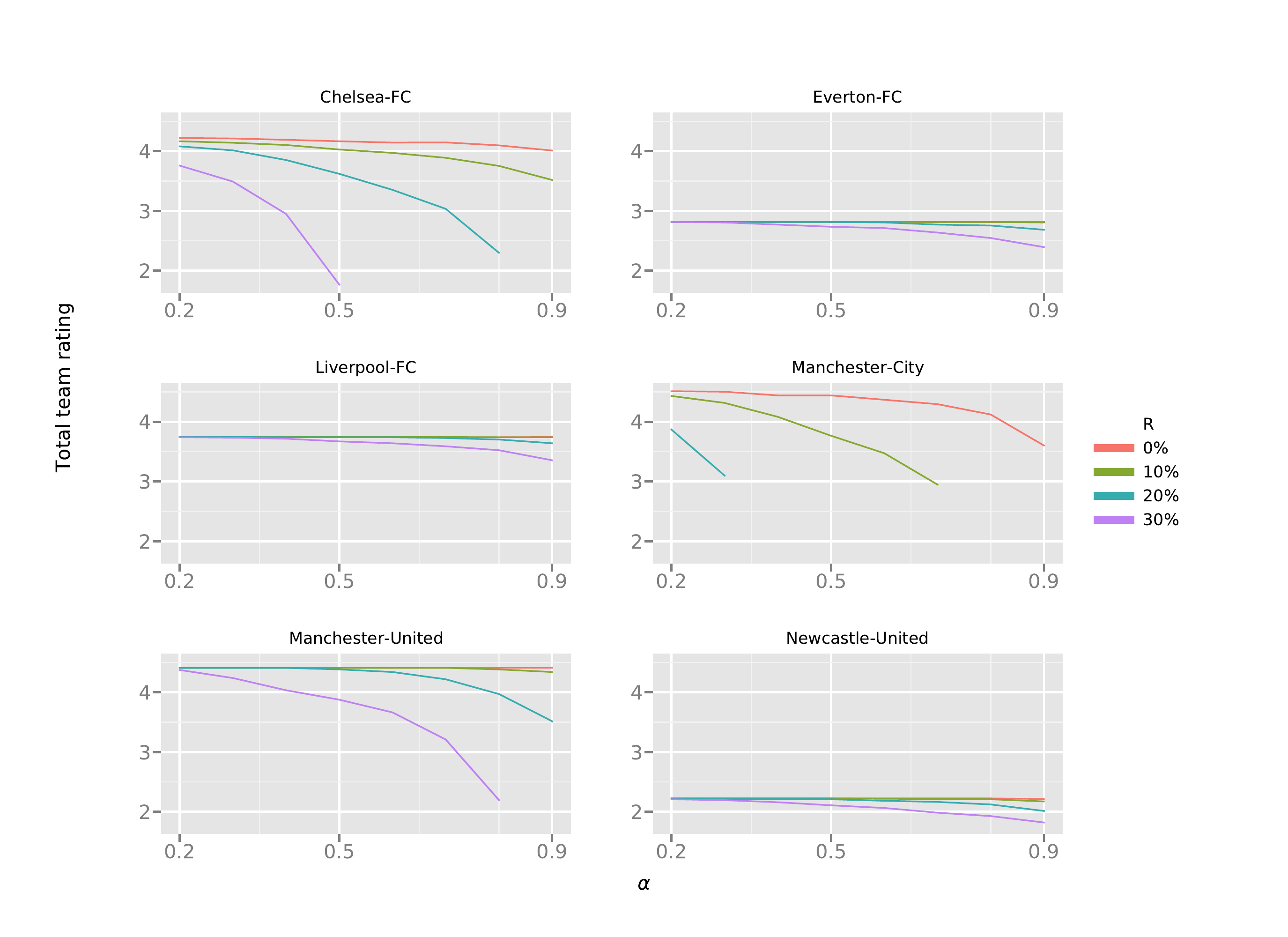}
\caption{Total rating of the resulting team for different degrees of risk tolerance $\alpha$ and for different values of $R$.}\label{fig:RatingVsAlpha}
\end{figure}

\Cref{fig:RatingVsAlpha} shows, for a sample of six clubs, an intuitive general trend: as the required probability of meeting financial goals increases,   the total team rating tends to decrease. As $\alpha$ increases we impose a higher probability of satisfying a purely financial measure. Consequently, the club has less freedom to sign top-performers, and is bound to find players that ensure a sufficient growth of the team value. Small $\alpha$ values represent clubs that are primarily interested in the here-and-now performance, and less concerned about the financial aspects. In this case, the decision maker has more freedom to choose top performers. 

For a number of clubs, pursuing a team value growth is incompatible with ensuring top-performers to the team, see, e.g., the case of Manchester City in \Cref{fig:RatingVsAlpha} with $R$ greater than $10\%$). However, a few clubs are rather insensitive to $\alpha$, especially for low values of $R$. 
In the latter case, the players that ensure the highest rating are, in general, the same player that ensure financial goals are met with sufficiently high probability. This is indeed a favorable situation, and it depends on the initial composition of the team as well as on the list of targets, and thus on the players available on the market. Notice, for example how Chelsea FC, Manchester City, and Manchester United show a similar high sensitivity to $\alpha$ as they share, in our case studies, the same list of target players. The same applies to Liverpool FC, Newcastle United, and Everton FC.

\begin{figure}[h!]
    \centering
    \begin{subfigure}{0.8\textwidth}
        \includegraphics[width=\textwidth]{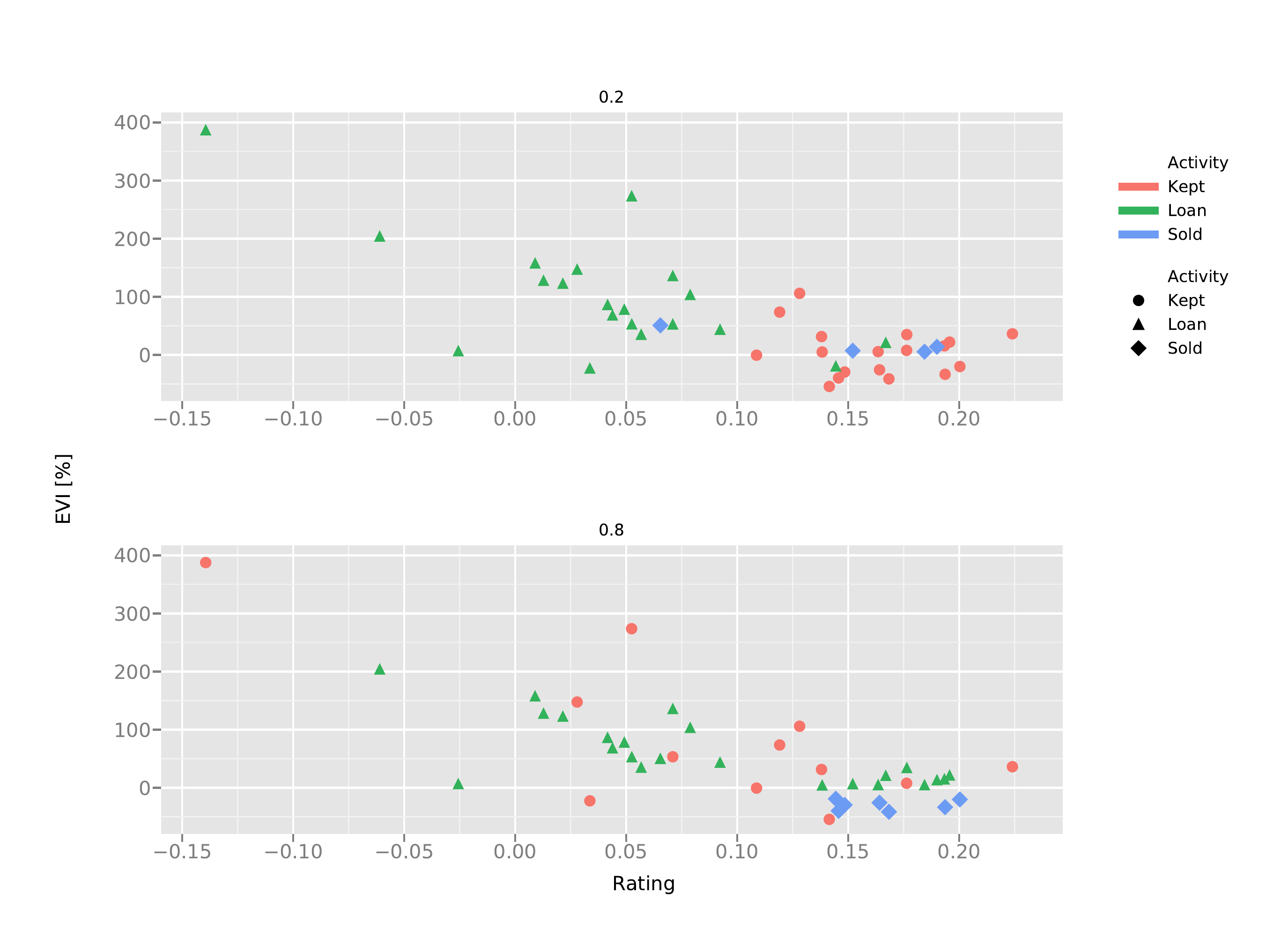}
        \caption{Transfers of own players}
        \label{fig:own02}
    \end{subfigure}\\
    \begin{subfigure}{0.8\textwidth}
        \includegraphics[width=\textwidth]{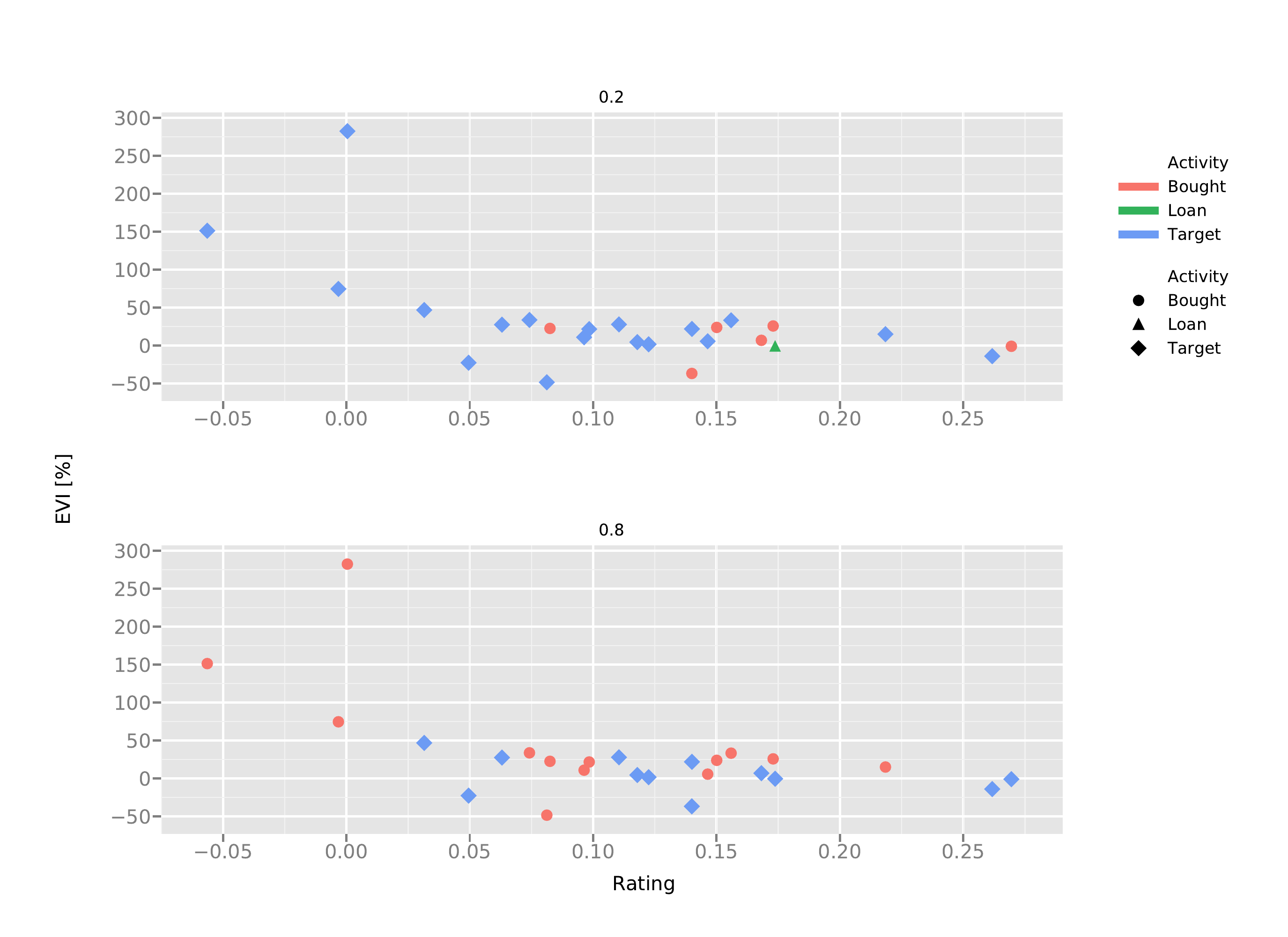}
        \caption{Transfers of target players}
        \label{fig:tar02}
    \end{subfigure}
    \caption{Distribution of transfers of Chelsea FC  $\alpha=0.2$ and $0.8$, and $R=1.2$. EVI represents the Expected Value Increase of a player, and is calculated as the expectation of \texttt{(future value-current value)/current value}.}\label{fig:chelsea}
\end{figure}

Let us zoom in on the case of Chelsea FC as a representative case for the clubs that are most sensitive to the probability $\alpha$. \Cref{fig:chelsea} reports the rating and the expected market value increase for the suggested transfers with $\alpha=0.2$ and $\alpha=0.8$ (assuming $R=1.2$). With $\alpha=0.2$ the club keeps most of the high-rating players, loans out most of the player with high expected growth that are not yet top performers, and sells some top performers with low expected value increase. On the other hand, when $\alpha=0.8$ the club will keep more of the high-expected-growth players and sell most of the high-rating players with low expected growth. Regarding inbound transfers, when $\alpha=0.2$ the club tends to buy, or sign on a loan agreement, players with above-average rating, and relatively low expected growth. However, with $\alpha=0.8$ the club signs the players with the highest expected value increase and fewer players with high ratings. That is, as the club becomes more concerned with financial stability, it will tend to build a team of high-potential players in spite of a reduced here-and-now performance. However, when the club is less concerned about finances, it will tend to keep its top performers and sign new high-rating players, in spite of the limited expected market value growth. 

\begin{figure}[h!]
    \centering
    \includegraphics[width=\textwidth]{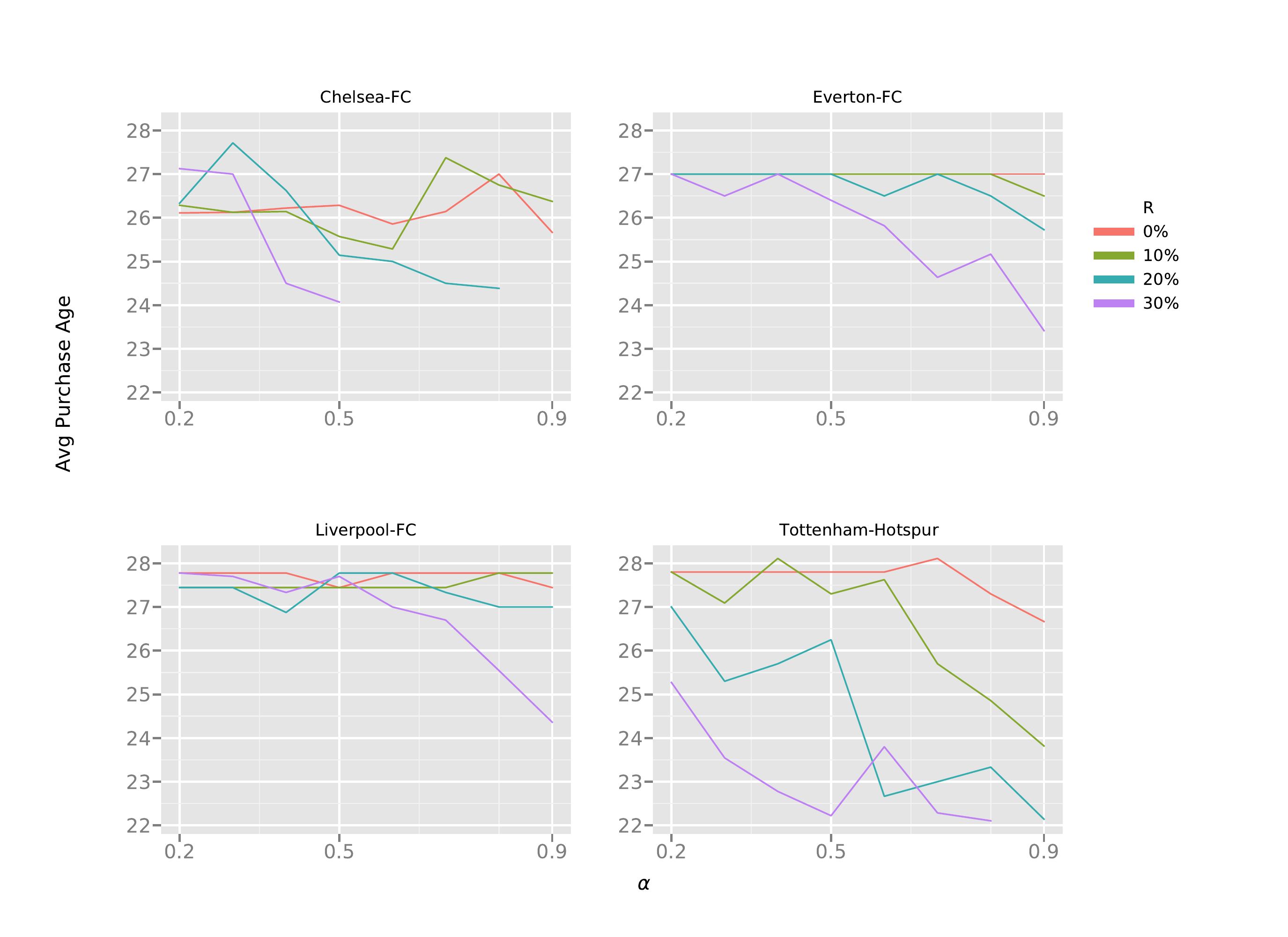}
    \caption{Average purchase age for different values of $\alpha$ and $R$.}\label{fig:purchaseVsAlpha}
\end{figure}
\clearpage
\begin{figure}[h!]
    \centering
    \includegraphics[width=\textwidth]{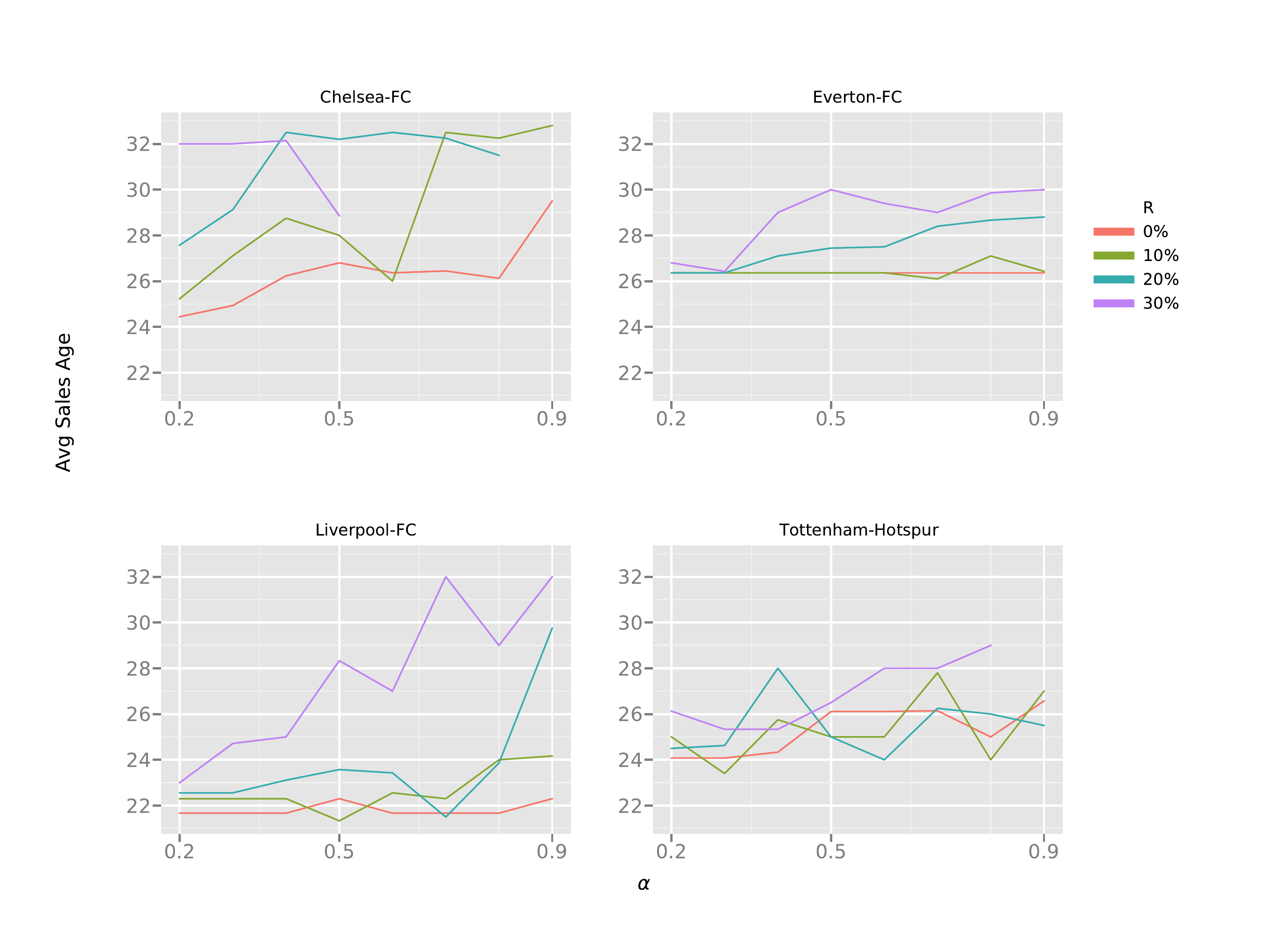}
    \caption{Average sales age for different values of $\alpha$ and $R$.}\label{fig:salesVsAlpha}
\end{figure}
\clearpage
As the requirement of meeting financial goals becomes stricter, the model suggests that Chelsea FC should buy younger players and sell older players, as illustrated in \Cref{tab:cs:chelseaAge}. The market value of younger players is, in general, expected to grow more than that of older players. Consequently, the model will tend to discard players that do not contribute to fulfilling the financial constraints. The insight we obtain from the solution for Chelsea FC is consistent with a more general trend that sees the average purchase age decrease with $\alpha$, and the average sales age increase with $\alpha$. These trends, illustrated in \Cref{fig:purchaseVsAlpha} and \Cref{fig:salesVsAlpha} for a sample of four clubs, are however erratic as purchases and sales need to generate a team composition feasible with respect not only to the chance constraints, but also to the other constraints of the problem. 

\begin{table}
  \centering
  \caption{Average age for the solution of Chelsea FC for $\alpha=0.2$ and $0.8$ and $R=1.2$. }
  \label{tab:cs:chelseaAge}
  \begin{tabular}{ccccc}
    \toprule
    $\alpha$ & Purchases & Hired on a loan agreement & Sales & Sent on a loan agreement \\
    \midrule
    0.2 & 26.0 & 28.0 &28.2 &20.6\\
    0.8 & 24.3 & --      &30.7    &23.0\\
    \bottomrule
  \end{tabular}
\end{table}

\subsection{Loaning strategies}\label{sec:results:loans}

\begin{figure}[htb]
    \centering
    \includegraphics[width=\textwidth]{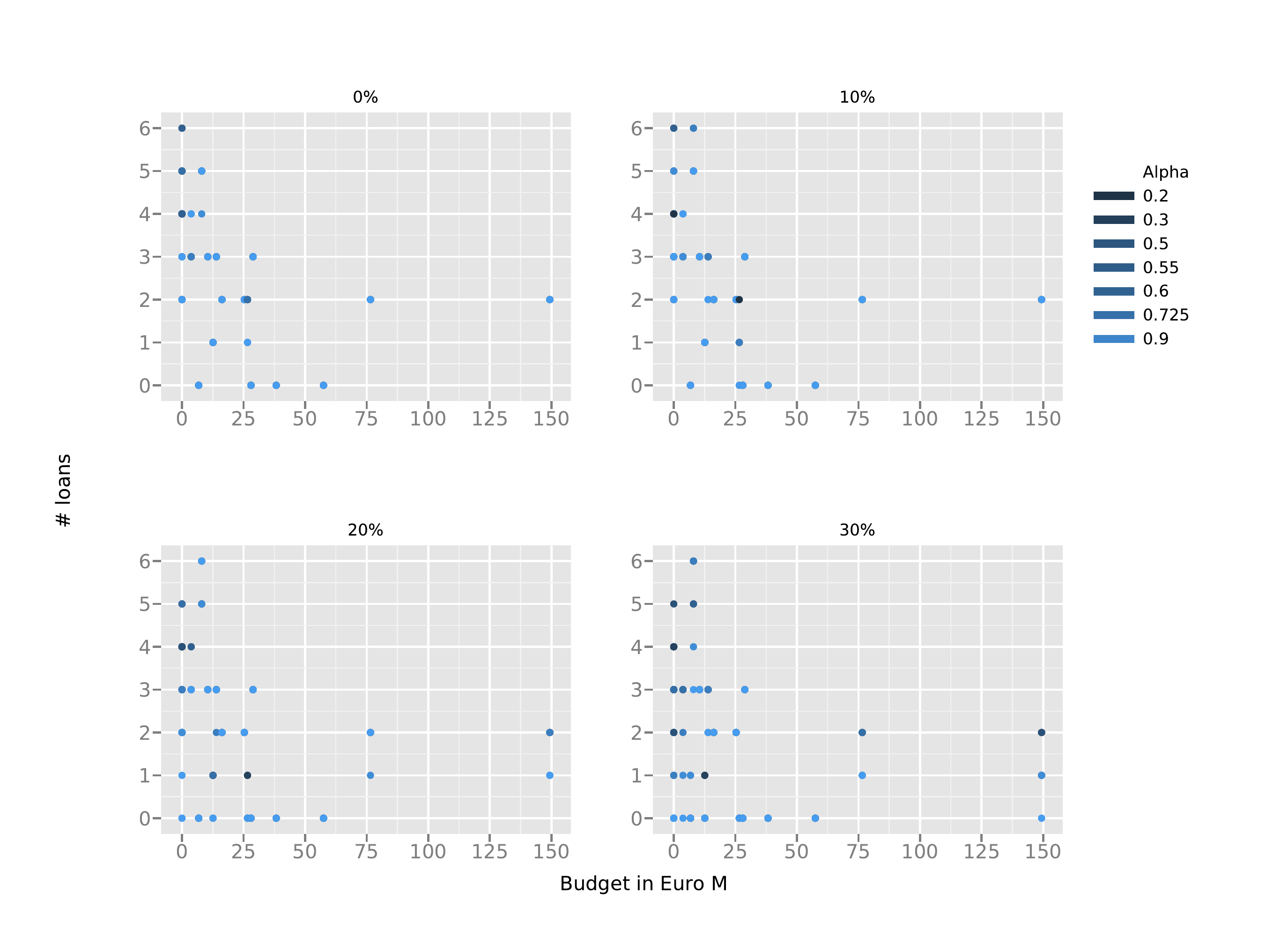}
    \caption{Number of players hired on a loan agreement for different budgets and values of $R$.}\label{fig:budget-loan}
\end{figure}

Hiring players on a loan agreement is a typical strategy for mid- and low-tier clubs to ensure a team of acceptable quality with a low budget to spend on the market. On the other hand, top clubs tend to purchase the players they need, very often due to more generous budgets. As shown in \Cref{fig:budget-loan}, the results obtained with model \eqref{eq:FTCP} are consistent with this general trend. The clubs that hire most players on a loan agreement are those with smaller budgets.

\subsection{Problem size and complexity}\label{sec:results:size}
Finally, we report on the size and complexity of the resulting optimization models. 
The size of the problems we solved, which is determined by the size of the sets $\mc{P}$, $\mc{R}$ and $\mc{S}$, is reported in \Cref{tab:size}.

\begin{table}[h]

    \centering
    \caption{Size of the problems. All variables are binary.}
    \label{tab:size}
    \begin{tabular}{c|cc}
    \toprule
        Team & \# Variables & \# Constraints  \\
        \midrule
         Arsenal-FC&275&272\\
Aston-Villa&345&329\\
Cardiff-City&375&348\\
Chelsea-FC&335&308\\
Crystal-Palace&370&344\\
Everton-FC&295&293\\
Fulham-FC&325&314\\
Hull-City&345&331\\
Liverpool-FC&300&296\\
Manchester-City&255&261\\
Manchester-United&280&275\\
Newcastle-United&320&308\\
Norwich-City&320&315\\
Southampton-FC&320&314\\
Stoke-City&320&312\\
Sunderland-AFC&325&314\\
Swansea-City&345&330\\
Tottenham-Hotspur&310&302\\
West-Bromwich-Albion&325&319\\
West-Ham-United&325&315 \\
\bottomrule
    \end{tabular}
\end{table}

All problems models were solved using the Java libraries of Cplex 12.6.2 on a machine equipped with CPU 2x2.4GHz AMD Opteron 2431 6 core and 24 Gb RAM. All problems were solved to a target 0.5\% optimality gap (parameter \texttt{EpGap} $0.5/100$) and using Cplex's default $10^{-6}$ feasibility tolerance and $10^{-5}$ integrality tolerance. These tolerances assume that values below $10^{-6}$ or $10^{-5}$ are either caused by rounding errors, or do not significantly change the optimization results. We set a time limit of $3600$ seconds (parameter \texttt{TimeLimit} $3600$) and put emphasis on proving optimality (parameter \texttt{MIPEmphasis} $3$). Descriptive statistics about the solution time in all our tests are reported in \Cref{tab:solution_time}. It can be noticed that the solution time is relatively contained for the great majority of the instances, offering space for further enhancements to the model.

\begin{table}[h]

    \centering
    \caption{Solution time statistics for the CC model across all the tests performed (640 model runs).}
    \label{tab:solution_time}
    \begin{tabular}{cc}
    \toprule
    Statistic & Value\\
    \midrule
    Average&        8.21 sec.\\
    St. dev.&       16.16 sec.\\
    Minimum &        0.04 sec.\\
    25th percentile &        0.25 sec.\\
    50th percentile &       1.32 sec.\\
    75th percentile &      10.53 sec.\\
    Maximum &      224.72 sec.\\
    \bottomrule
\end{tabular}
\end{table}

\section{Concluding remarks}
\label{sec:Conclusions}
This article introduced a chance-constrained optimization model for assisting football clubs during transfer market decision. Furthermore it presented a new rating system which is able to measure numerically the on-field performance of football players. Such measure is necessary in order to arrive at an objective assessment of football players and thus limit the bias in the observers. 

The model and rating system have been extensively tested on case-studies based on real-life English Premier League marked data. The results illustrate that the model  contribute to reduce bias in transfer market decision by supporting football managers' expertise with analytic support and validation. Furthermore, the model can adapt to different levels of financial concern and thus support football decision makers with tailor-made analytic suggestions. 

There is still room for improvements and enhancements of the model and the rating system. For example, the rating system only considers the performance of players when appearing on the pitch. However, some players may be prone to injuries or suspensions, so that their use in a team may be limited.

\subsection*{Acknowledgements}

The authors wish to thank two anonymous reviewers, whose comments helped to improve the contents and presentation of this manuscript.


\bibliography{bibfile}

\appendix

\section{Sensitivity to the formation}\label{app:formation}
\begin{figure}[htb]
    \centering
    \includegraphics[width=\textwidth]{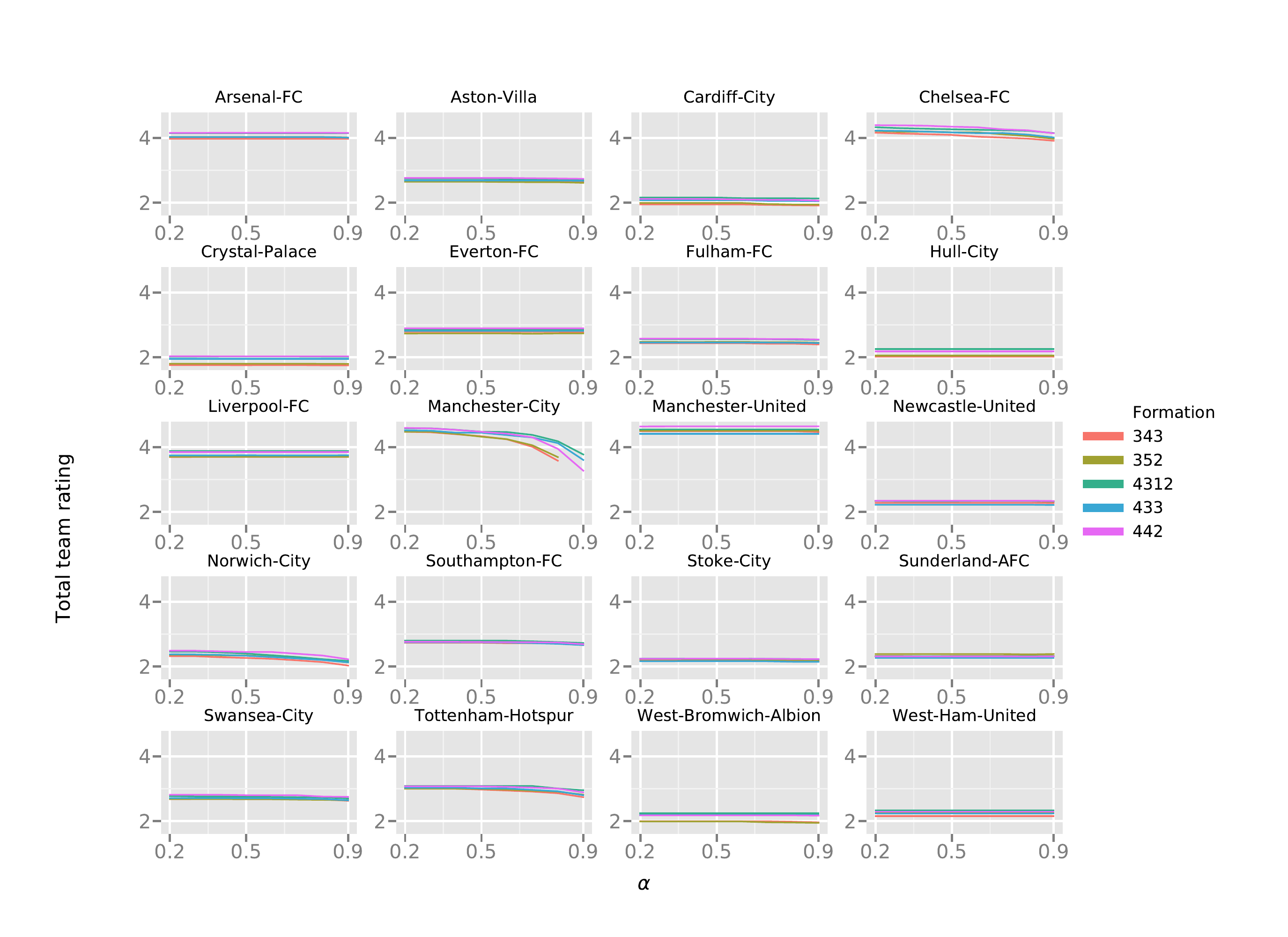}
    \caption{Total team rating for different values of $\alpha$ and for different formations with $R=1.0$.}\label{fig:formation}
\end{figure}

\Cref{fig:formation} illustrates that solutions to model \eqref{eq:FTCP} are not sensitive to the formation, and thus the discussion in \Cref{sec:results} holds independently of the formation chosen.

\end{document}